\documentclass[12pt]{amsart}
\usepackage{amssymb,latexsym}
\usepackage[mathscr]{eucal}
\usepackage{graphics}
\usepackage{ulem}
\usepackage{verbatim}

\theoremstyle{plain}
\newtheorem{theorem}{Theorem}[section]
\newtheorem{corollary}[theorem]{Corollary}
\newtheorem{lemma}[theorem]{Lemma}
\newtheorem{proposition}[theorem]{Proposition}
\newtheorem{remark}[theorem]{Remark}
\newtheorem{definition}[theorem]{Definition}
\newtheorem{example}[theorem]{Example}
\newtheorem{question}[theorem]{Question}
\DeclareMathOperator{\diag}{diag}
\DeclareMathOperator{\und}{und}
\DeclareMathOperator{\lnd}{lnd}
\DeclareMathOperator{\uni}{uni}
\DeclareMathOperator{\lni}{lni}
\DeclareMathOperator{\se}{se}
\DeclareMathOperator{\scop}{sc}
\DeclareMathOperator{\cl}{cl}
\DeclareMathOperator{\di}{dim}

\def \co*{{\text{c}_{\text{o}}^*}}
\def \D*{$\Delta_{1/2}$-condition}
\newcommand{\be}{\begin{equation}\label}
\newcommand{\ee}{\end{equation}}
\newcommand{\bq}{\begin{equation*}}
\newcommand{\eq}{\end{equation*}}
\newcommand{\ba}{\begin{align*}}
\newcommand{\ea}{\end{align*}}
\newcommand{\bp}{\begin{proof}}
\newcommand{\ep}{\end{proof}}
\newcommand{\bL}{\begin{lemma}\label}
\newcommand{\eL}{\end{lemma}}
\newcommand{\bP}{\begin{proposition}\label}
\newcommand{\eP}{\end{proposition}}
\newcommand{\bC}{\begin{corollary}\label}
\newcommand{\eC}{\end{corollary}}
\newcommand{\bT}{\begin{theorem}\label}
\newcommand{\eT}{\end{theorem}}
\newcommand{\bR}{\begin{remark}\label}
\newcommand{\eR}{\end{remark}}
\newcommand{\bD}{\begin{definition}\label}
\newcommand{\eD}{\end{definition}}
\newcommand{\bQ}{\begin{question}\label}
\newcommand{\eQ}{\end{question}}

\begin{document}

\title[Arithmetic mean ideals]
{A survey on the interplay between \quad \\
arithmetic mean ideals, traces, \quad \\
lattices of operator ideals, and \quad \\ 
an infinite Schur-Horn majorization theorem}
\author{Victor Kaftal}
\address{University of Cincinnati\\
          Department of Mathematics\\
          Cincinnati, OH, 45221-0025\\
          USA}
\email{victor.kaftal@math.uc.edu}
\author{Gary Weiss}
\email{gary.weiss@math.uc.edu}

\date{April 3, 2007}

\begin{abstract}

    The main result in \cite{DFWW} on the structure of commutators showed that arithmetic means play an important role in the study of operator ideals. In this survey we present the notions of arithmetic mean ideals and arithmetic mean at infinity ideals. Then we explore their connections with commutator spaces, traces, elementary operators, lattice and sublattice structure of ideals, arithmetic mean ideal cancellation properties of first and second order, and softness properties - a term that we introduced but a notion ubiquitous in the literature on operator ideals. Arithmetic mean closure of ideals leads us to investigate majorization for infinite sequences and this in turn leads us to an infinite Schur-Horn majorization theorem
which extends theorems by A. Neumann, by Arveson and Kadison, and by Antezana, Massey, Ruiz and Stojanoff. This survey covers the material announced towards the beginning of the project in PNAS-US \cite{vKgW02} and then expanded and developed in a series of papers  \cite{vKgW04-Traces}-\cite{vKgW04-Majorization}. We also list ten open questions that we encountered in the development of this material.
\end{abstract}

\maketitle

\section{\leftline{\bf Introduction}}\label{S:1}

Commutator spaces (also known as commutator ideals) of operator ideals were characterized in \cite{DFWW} in terms of arithmetic mean operations and this led already in \cite{DFWW} to the notion of two basic arithmetic mean ideals, the arithmetic mean closure and the arithmetic mean interior of an ideal, and to the notion of arithmetically mean stable ideals - those that coincide with their commutator space and hence support no nonzero trace. Thus  arithmetic mean ideals have become an intrinsic part of the theory of operator ideals. To explore their structure and to find how they relate to the existing theory of operator ideals and in particular, how they contribute to it, was the goal of our program.  

The first results of this work were announced in PNAS-US \cite{vKgW02} and then expanded and developed in a series of papers \cite{vKgW04-Traces}-\cite{vKgW04-Majorization}. 

Our beginning question, explored in \cite{vKgW04-Traces}, was to ask ``how many traces can a given ideal support."  To address it, we found that we first needed both to investigate much more systematically the structure of arithmetic mean ideals and to introduce formally the notion of ``soft ideals." 

From the start we discovered that, to study the two basic arithmetic mean ideals, we had to expand our investigations by adding to the menagerie the smallest am-open ideal containing an ideal and the largest am-closed ideal contained by it. The latter, for instance, provided a useful tool for the study of countably generated am-stable ideals. 

The arithmetic mean is too crude for distinguishing ideals contained in the trace class $\mathscr L_1$: for them what is required is the arithmetic mean at infinity. 
The arithmetic mean at infinity (am-$\infty$ for short) has substantially different properties from the arithmetic mean, e.g., the am-$\infty$ of a sequence can fail to satisfy the $\Delta_{1/2}$-condition whereas an important feature of the arithmetic mean of a sequence is that it always satisfies this condition.  
The ensuing theory of am-$\infty$ ideals is considerably more complicated, but in large measure mirrors the arithmetic mean case. This material is presented in Section \ref {S:2} and was developed in  \cite{vKgW04-Traces} and \cite{vKgW04-Soft}.

In Section \ref{S:3} we present ``soft ideals" which are ubiquitous in the literature on operator ideals and which we studied systematically in 
\cite{vKgW04-Soft}. We show that many of the classical ideals are soft and that keeping this in mind can add perspective and simplify the proofs of a number of results in the literature.  
Our interest in the notion of softness arose from the interplay between softness properties and arithmetic mean operations, which provided the main tool for our work on traces. We explore in this section this interplay and present several open questions.

Section \ref{S:4} is devoted to traces and based is on \cite{vKgW04-Traces}.
The question "how many traces can a given ideal support"  asks about the codimension of the commutator space $[I,B(H)]$ of an ideal $I$. 
We were led by our work in    \cite{vKgW04-Traces} to the conjecture that
\[
\di \frac{I}{[I,B(H)]} ~\in
\begin{cases}
\{1,\infty\}    &\text{when}~\omega \notin \Sigma(I) \\
\{0,\infty\}    &\text{when}~\omega \in \Sigma(I)
\end{cases}
\]
where $\omega:= \,<1, \frac{1}{2}, \frac{1}{3}, \dots>$ is the harmonic sequence and $\Sigma(I)$ is the characteristic set of $I$. 
In the case when $\omega \notin \Sigma(I)$, we show that the codimension is $1$  precisely when $I$ is am-$\infty$ stable. In the process, we obtain results that can be applied to solve questions on elementary operators. \\

\noindent All ideals divide naturally into the classes:
\begin{itemize}
\item  ``small ideals," 
(those contained in the lower stabilizer at infinity of the principal ideal $(\omega)$, $st_{a_\infty}(\omega)$, i.e., 
the ideal with characteristic set $\{\xi \in \co* \mid \sum \xi_nlog^m n < \infty\}$), 
\item ``large ideals" (those that contain the upper stabilizer $st^a(\omega) :=\bigcup^\infty_{m=0}(\omega\log^m)$, a directed union of the principal ideals $(\omega\log^m)$)  and 
\item``intermediate ideals" (all the remaining ideals)

\end{itemize}

\noindent For all intermediate ideals, the codimension of the commutator space  is always infinite. Then we show that the above codimension conjecture holds for all soft ideals and present other sufficient conditions for the codimension of the commutator space to be infinite.

Section \ref{S: Density} is based on \cite{vKgW04-Density}, a study of the lattice structure of operator ideals focusing in particular on some distinguished sublattices of principal ideals. We investigate density properties, i.e., when between two ideals each in one of two nested lattices lies another from the smaller lattice. These properties are powerful tools for dealing with general ideals. We also study representations of ideals as directed unions of countably generated ideals and of principal ideals.

We apply then some of these tools to problems on arithmetic mean ``equality cancellations": 
when for an ideal $I$ does $J_a=I_a\Rightarrow J=I$? 
And ``inclusion cancellations": when for an ideal $I$ does $J_a\subset I_a\Rightarrow J\subset I$ or $J_a\supset I_a\Rightarrow J\supset I$? 
The most interesting case is the latter and the answer is when $I=\widehat{I}$. 
For $I$ principal and not contained in the trace class we show that $\widehat{I}$ is principal as well and we prove that for for $0<p<1$, 
$\widehat{(\omega^p)}=(\omega^{p'})$  where $\frac{1}{p} -\frac{1}{p'}=1$. 
We show that the condition $(\xi)= \widehat{(\xi)}$ is strictly stronger than regularity.

Second order cancellation properties are considerably different and harder than first order cancellation properties. 
We found that the cancellation $I_{a^2}=J_{a^2}\Rightarrow J_a=I_a$ does not hold in general even if $I$ and $J$ are principal ideals, 
thus answering a question by M. Wodzicki. 
Sufficient conditions for this and for the two second order inclusion cancellations to hold are given in Section \ref{S: 2nd Order} in terms of the "ratio of regularity," 
$r(\xi_a): = \frac{\xi_{a^2}}{\xi_a}$, where $\xi$ is the generator of $I$. This material was developed in \cite {vKgW04-2nd Order}.

Finally, Section \ref{S: Majorization} focuses on majorization theory for infinite sequences. The link with operator ideals is that the definition of am-closure of an ideal can be restated in terms of hereditariness (i.e., solidity) with respect to majorization.  
Our initial goal was to prove that am-closure for an ideal is equivalent to diagonal invariance, i.e., the property that for any fixed orthonormal basis, if an operator $A$ belongs to the ideal then its main diagonal, $E(A)$, also belongs to the ideal. 
To prove this equivalence, we introduced the notions of block majorization and strong majorization of sequences and these, in turn, provided the tools and the stimulus for proving an infinite Schur-Horn majorization theorem (Theorem \ref{T:New Horn}). 
The key part of this theorem is that if $\xi$ and $\eta$ are monotone sequences decreasing to $0$ and $\xi$ is not summable, then the condition 
$\sum_{j=1}^n\xi_j\le \sum_{j=1}^n\eta_j$ for every $n$ is necessary and sufficient for $\diag \xi = E(U\diag\eta U^*)$ for some unitary operator $U$. 
The material of  Section \ref{S: Majorization} is based on \cite {vKgW04-Majorization} which is not yet in final form but we expect it to be available by the time this survey is published. \\

No proofs are presented  in this survey and many of the basic facts are stated without explicit reference, but both  proofs and further background information can be found easily in the respective citations. \\

We include herein ten questions that arose in this project.

\section{\leftline{\bf Commutators, traces and arithmetic mean ideals}}\label{S:2}
The natural domain of the usual trace $Tr$ on $B(H)$ (with $H$ a separable infinite-dimensional complex Hilbert space) 
is the trace class ideal $\mathscr{L}_1$, but ideals of $B(H)$ can support other traces.

\begin{definition}\label{D:trace} 
A trace $\tau$ on an ideal $I$ is a unitarily invariant linear functional on $I$. 
\end{definition}
\noindent In this paper, traces are neither assumed to be positive nor faithful. (We refer the reader to \cite {jV89} and the more recent  \cite{mW02} for results on positive traces.) \textit{All ideals are assumed to be proper and two-sided.} 
All proper ideals contain the finite rank ideal $F$ and the restriction of any trace to $F$ is a scalar multiple $cTr$ of the standard trace $Tr$. 
\textit{Singular traces} are those for which $c=0$ and \textit{nonsingular traces} are those for which $c\ne 0$. 

Recall that commutators are operators $[X,Y]:= XY-YX$ and that $[I,B(H)]$ denoting the \textit{commutator space} 
(also known as the \textit{commutator ideal}) of $I$ and $B(H)$ is defined as the linear span of commutators, $XY-YX$, of operators $X$ in $I$ with operators $Y$ in $B(H)$ (and likewise for $[I,J]$). 

Since $UXU^* - X = [U, XU^* - U^*X] \in [I,B(H)]$ for every $X \in I$ and every unitary operator $U$ and since unitary operators span $B(H)$,
unitarily invariant linear functionals on an ideal $I$ are precisely the linear functionals on $I$ that vanish on $[I,B(H)]$. 
Thus traces can be identified with the elements of the linear dual of the quotient space $\frac{I}{[I,B(H)]}$, 
which is why commutator spaces play such a central role in the study of traces.
The structure of general $[I, J]$ had been studied for several decades \cite{pH54}, \cite{BP65}, \cite{PT71}, \cite{jA77}, \cite{gW75}-\cite{gW86} and \cite{nK89} (see also \cite{gW05} for other historical references).  
The introduction of cyclic cohomology in the early 1980's by 
A. Connes and its linkage to algebraic K-theory by 
M. Wodzicki in the 1990's 
provided additional motivation for the complete determination of commutator ideal structure.
(Cf. \cite{aC82}-\cite{aC85}  and \cite{mW94}.) \\

Clearly the commutator spaces of ideals are selfadjoint linear spaces and hence can be characterized via their selfadjoint elements and those (and more generally the normal operators in a commutator space) were  characterized in \cite{DFWW} in terms of arithmetic means. 
An important feature is that membership in commutator spaces (noncommutative objects) 
is reduced to certain conditions on associated sequences and their arithmetic means (commutative objects).

Arithmetic means first entered the analysis of $[\mathscr L_2,\mathscr L_2]$ and $[\mathscr L_1,B(H)]$ 
for a special case in \cite{gW75} and \cite{gW80} and for 
its full characterization in \cite{nK89}.  
As the main result in \cite{DFWW} (Theorem 5.6) conclusively shows, 
arithmetic means are essential for the study of traces and commutator spaces in operator ideals. 
\cite{DFWW} also initiated a systematic study of ideals derived via the arithmetic mean operations and in \cite{vKgW02}-\cite{vKgW04-2nd Order}
we continued this study along with the arithmetic mean at infinity operations and the consequent parallel theory described below.\\

When $X \in K(H)$ (the ideal of compact operators on $H$), denote an ordered spectral sequence for $X$ by $\lambda(X) :=~ <\lambda(X)_1,\lambda(X)_2,\dots>$, 
i.e., a sequence of all the eigenvalues of $X$ (if any), repeated according to algebraic multiplicity, 
completed by adding infinitely many zeroes when only finitely many eigenvalues are nonzero, and arranged in any order so that $|\lambda(X)|$ is nonincreasing.  
For any sequence $\lambda =~ <\lambda_n>$, denote by $\lambda_a$, $\lambda_{a_\infty}$ the sequences of its arithmetic (Cesaro) mean and arithmetic mean at infinity, i.e.,
\[
\lambda_a :=~ \left <\frac{1}{n}\sum_{j=1}^{n}\lambda_j \right>_{n=1}^{\infty} \quad \text{, and when}~ \lambda \in \ell^1, \quad
\lambda_{a_\infty} :=~ \left <\frac{1}{n}\sum_{j=n+1}^{\infty}\lambda_j \right>_{n=1}^{\infty}. 
\]

A constant theme in the theory of operator ideals 
is its connection to the theory of sequence spaces.  
Calkin \cite{jC41} established a correspondence between the two-sided ideals of $B(H)$
and the \textit{characteristic sets}, i.e., the positive cones of $\text{c}_{\text{o}}^*$ 
(the collection of sequences decreasing to $0$)
that are hereditary (solid under pointwise order) and invariant under ampliations 
\[
\text{c}_{\text{o}}^* \owns \xi \rightarrow D_m\xi:=~<\xi_1,\dots,\xi_1,\xi_2,\dots,\xi_2,\xi_3,\dots,\xi_3,\dots>
\] 
where each entry $\xi_i$ of $\xi$ is repeated $m$-times. 
The order-preserving lattice isomorphism 
$I \rightarrow \Sigma(I)$ maps each ideal to its characteristic set $\Sigma(I) := \{s(X) \mid X \in I\}$ where $s(X)$ denotes the sequence of 
$s$-numbers of $X$, i.e., all the eigenvalues of $|X| = (X^*X)^{1/2}$ repeated according to multiplicity, arranged in decreasing order, and completed by adding infinitely many zeroes if $X$ has finite rank.  
Conversely, for every characteristic set $\Sigma \subset \text{c}_{\text{o}}^*$, 
if $I$ is the ideal generated by $\{\diag \xi \mid \xi \in \Sigma\}$ 
where $\diag \xi$ is the 
diagonal matrix with entries $\xi_1,\xi_2,\dots$, then we have $\Sigma = \Sigma(I)$.  

For ideals $I,J$, the ideals $I+J$, $IJ=JI$ and $I:J$ are, respectively, the ideal of sums of operators from $I$ and $J$, 
the ideal generated by all products of operators from $I$ and $J$ (single products will do \cite[Lemma 6.3]{DFWW}), 
and the ideal of all operators whose products with all operators in $J$ are in $I$.
In terms of characteristic sets, 
\[
\Sigma(I+J)=\Sigma(I)+\Sigma(J)\quad \text{(from } s_{2n+1}(A+B) \le s_n(A)+s_n(B)\text{ \cite{GK69}), }
\]
\[
\Sigma(IJ)=\{\xi\in \co*\mid \xi \le \eta \zeta\,\text{ for some}\,\eta \in \Sigma(I), \zeta\in \Sigma(J)\}
\] 
(hence $IJ$ is a commutative operation), and 
\[
\Sigma(I:J)=\{\xi\in \co*\mid \xi\,\Sigma(J) \subset \Sigma(I)\}. 
\]

A special case of \cite[Theorem 5.6]{DFWW} (see also [ibid, Introduction]) characterizing all commutator ideals $[I,J]$ 
and yielding the important identity $[I,J]=[IJ,B(H)]$, 
is the characterization of the normal elements of the commutator space $[I, B(H)]$: 

\begin{theorem} \label{T:DFWW}\cite [Theorem 5.6]{DFWW} 
Let $X \in I$ be a normal operator in a proper ideal $I$ with any ordered spectral sequence $\lambda(X)$.
Then $X \in [I,B(H)]$ if and only if $|\lambda(X)_a| \leq \xi$ for some $\xi \in \Sigma(I)$.
\end{theorem}

Two important consequences are: \\
(i) The introduction in \cite{DFWW} of arithmetic mean ideals to characterize $[I, B(H)]^+$, the positive part of  $[I, B(H)]$.\\
(ii) The standard trace $Tr$ extends from $F$ to $I$ if and only if $\diag \omega \notin I$ where the harmonic sequence $\omega := \,<1,\frac{1}{2},\frac{1}{3},\dots>$.

The action of arithmetic means on sequences indeed leads naturally to an action on operator ideals.
The four basic arithmetic mean ideals called respectively the pre-arithmetic mean, arithmetic mean, pre-arithmetic mean at infinity, and arithmetic mean at infinity of $I$ are $_aI$, $I_a$, $_{a_\infty}I$, $I_{a_\infty}$ (see \cite[Section 2.8(v)-(vi)]{DFWW} and \cite[Definition 4.6]{vKgW04-Traces}),  and their characteristic sets are defined as:
\[
\Sigma(_aI) := \{\xi \in \text{c}_{\text{o}}^* \mid \xi_a \in \Sigma(I)\} 
\]
\[
\Sigma(I_a) := \{\xi \in \text{c}_{\text{o}}^* \mid \xi = O(\eta_a)~\text{for some}~ \eta \in \Sigma(I)\}
\]
\[
\Sigma(_{a_\infty}I) := \{\xi \in (\ell^1)^* \mid \xi_{a_\infty} \in \Sigma(I)\} \quad \text{(defined only for $I \ne \{0\}$)}
\]
\[
\Sigma(I_{a_\infty}) := \{\xi \in \text{c}_{\text{o}}^{*} \mid \xi = O(\eta_{a_\infty})  ~\text{for some}~ \eta \in \Sigma(I \cap \mathscr L_1)\}.
\]
Notice that in this paper  $\co*$ and $(\ell^1)^*$ denote the positive cones of monotone nonincreasing sequences in c$_o$ and $(\ell^1)^+$ and not the duals of those spaces. A special case of the main result of \cite{DFWW}, Theorem \ref{T:DFWW} above, restricted to positive operators, can be reformulated in terms of arithmetic mean ideals and another consequence is the am-$\infty$ analog \cite[Corollary 6.2(i)]{vKgW04-Traces}:

\begin{theorem} \label{T:pos part + infty characterization}
\item[(i)] For every proper ideal $I$: $[I,B(H)]^+ = (_aI)^+$
\item[(ii)]  For every ideal $I \subset \mathscr L_1$ (or more generally if $\omega \notin \Sigma(I)$): 
$(F+[I,B(H)])^+ = (_{a_\infty}I)^+$
\end{theorem}
An immediate and interesting consequence is that both $[I,B(H)]^+$  and $(F+[I,B(H)])^+$ are hereditary (solid relative to the order $A \le B$). 
Later we shall see that so also is $(\mathscr L_1+[I,B(H)])^+$ and ask about general $(J+[I,B(H)])^+$.

\subsection{\leftline{Stability, stability at infinity and stabilizers}}\label{s: stability and infty-stability}\quad
An important consequence of Theorems \ref{T:DFWW} is that $I=[I,B(H)]$ if and only if $I=\,_aI$. 
It is easy to show that $I=\,_aI$ if and only if $I=I_a$.
And similarly, though not exactly, for the am-$\infty$ case:

\bD{D:am stable} 
An ideal $I$ is called arithmetically mean stable (am-stable for short) if and only if $I = \,_aI$ if and only if $I = I_a$ if and only if $I=[I,B(H)]$.

An ideal $I$  is called arithmetically mean at infinity stable (am-$\infty$ stable for short) if and only if $I = \,_{a_\infty}I$ if and only if $\mathscr L_1 \supset I = I_{a_\infty}$ if and only if $\mathscr L_1 \supset I=F+[I,B(H)]$ 
\cite[Corollary 4.10 and Theorem 6.6]{vKgW04-Traces}. 
\eD

Am-stable ideals are precisely the ideals that do not support any nonzero trace, and over the years a considerable amount of effort has been invested in determining which ideals are stable: 
\cite{pH54}, \cite{GK69}, \cite{gW75}, \cite{gW80}, \cite{gW86}, \cite{AS78}, \cite{AV86}, \cite{jA86} and \cite{jV89} listed chronologically. 
Am-stability for many classical ideals was also studied in \cite[Sections 5.11--5.27]{DFWW}. \\

Of particular interest are principal ideals. 
If $A\in K(H)$, let $\xi:=s(A)$ be the s-number sequence of $A$, 
then $A$ and $\diag \xi $ generate the same principal ideal, which we will denote by $(\xi)$. 
It is easy to see that $(\xi)= (\eta)$ if and only if $\xi= O(D_m\eta)$ and $\eta = O(D_{m'}\xi)$ for some $m,m'\in \mathbb N$. 
$\xi$ is said to satisfy the \D* condition if $D_2\xi= O(\xi)$ (if and only if $D_m\xi= O(\xi)$ for every $m\in \mathbb N$). 
If $(\xi)= (\eta)$ and $\xi$ satisfies the \D* condition, then so does $\eta$ and $\xi \asymp \eta$. 
Sequences for which $\xi\asymp \xi_a$ are called ``regular" and have been studied in connection with principal ideals 
since Gohberg and Krein \cite[p.143 (14.12)]{GK69}. 
Reformulated in terms of am-stability therefore, a sequence $\xi$ is regular if and only if the principal ideal $(\xi)$ is am-stable.
Principal ideals that are am-$\infty$ stable are characterized in Theorem \ref {T:infty regular} below.

Every ideal $I$ is contained in the smallest am-stable ideal containing it, which we call the \textit{upper am-stabilizer} of $I$, 
$st^a(I) := \bigcup^\infty_{m=0}I_{a^m}$, and contains a largest (possibly zero) am-stable ideal which we call the \textit{lower am-stabilizer} of $I$, $st_a(I) :=\bigcap^\infty_{m=0}\,_{a^m}I$. 
It follows that there is no largest proper am-stable ideal and the smallest nonzero am-stable ideal is 
\[
 st^a(F) = st^a((\omega)) = st^a(\mathscr L_1) :=
\bigcup^\infty_{m=0}(\omega)_{a^m}=
\bigcup^\infty_{m=0}(\omega\log^m)  \supset \mathscr L_1
\]
\cite[Proposition~4.18]{vKgW04-Traces}. 

Similarly every ideal contains a largest am-$\infty$ stable ideal, which we call the \textit{lower am-$\infty$ stablilizer} of $I$, 
$st_{a_\infty}(I) := \cap^\infty_{m=0}\,_{a_\infty^m}(I)$.
Notice that every am-$\infty$ stable ideal must be contained in $\mathscr L_1$ and hence in
\[
st_{a_\infty}(K(H))=st_{a_\infty}((\omega))= st_{a_\infty}(\mathscr L_1) = \cap^\infty_{m=0}\,_{a_\infty^m}(\mathscr L_1) = \cap^\infty_{m=0}\,\mathscr L(\sigma\,log^m) \subset \mathscr L_1,
\]
where $\mathscr L(\sigma\,log^m)$ is the Lorentz ideal with characteristic set $\{\xi \in \co* \mid \sum \xi_nlog^m n < \infty\}$ 
\cite[Proposition 4.18]{vKgW04-Traces}.
For ideals $I$ contained in $st_{a_\infty}(\mathscr L_1)$, the smallest am-$\infty$ stable ideal containing $I$ is its \textit{upper am-$\infty$ stabilizer} $st^{a_\infty}(I) := \cup^\infty_{m=0}(I_{a_\infty^m})$.

The ideals $st^a(\mathscr L_1)$ and $st_{a_\infty}(\mathscr L_1) $ naturally divide all the ideals into three classes (see \cite[Introduction]{vKgW04-Traces}):
the ``small ideals," those contained in $st_{a_\infty}(\mathscr L_1)$, the ``large ideals," those containing $st^a(\mathscr L_1)$, 
and the ``intermediate ideals," those that are neither.
This provides a useful perspective from which to attack the codimension conjecture for traces (see Section \ref {S:2.2.2}).

Of special importance in \cite{DFWW} and herein is the principal ideal $(\omega)$. Elementary computations show that $F_a = (\mathscr L_1)_a = (\omega)$ and that $_a(\omega) = \mathscr L_1$. 
Hence $_aI \ne \{0\}$ if and only if $\omega \in \Sigma(I)$. 
An immediate but important consequence of Theorem \ref{T:DFWW} which is used often throughout \cite{vKgW04-Traces} is that 
$\omega \in \Sigma(I)$ if and only if $\mathscr L_1 \subset [I,B(H)]$ 
if and only if $F \subset [I,B(H)]$.\\

\subsection{\leftline{Am-closure and am-interior}}\label{s:clos}\quad
\cite{DFWW} defined for every ideal $I$ two further arithmetic mean ideals, 
the arithmetic-mean closure $I^-:=\,_a(I_a)$ and the arithmetic mean interior $I^o: = (_aI)_a$ which together with the ideals $I_a$ and $_aI$ we call ``first order arithmetic mean ideals". 
Both the am-closure and am-interior operations are idempotent and it is natural to call am-closed 
(resp., am-open) those ideals that coincide with their am-closure, i.e., $I=I^-$  (resp., am-interior, i.e., $I=I^o$). 
Thus $I^o$ is the largest am-open ideal inside $I$ and $I^-$ is the smallest am-closed ideal containing $I$.
(We do not mean to imply a full topological structure as, for instance, the union of am-open ideals can even fail to be an ideal.)
Many of the classical sequence space ideals in the literature are am-closed.
Am-closed ideals also play an important role in the study of single commutators \cite [Section 7]{DFWW}. 
The smallest am-closed ideal is the trace class $\mathscr L_1$ and the smallest  am-open ideal is $(\omega)$. 
There are no proper largest am-closed or am-open  ideals.

It is easy to see that an ideal $I$ is am-open (resp., am-closed) if and only if $I=J_a$ (resp., $I=\,_aJ$ for some ideal $J$ and that the following ``five chain of inclusions" holds:
\[
_aI\subset I^o\subset I\subset I^-\subset I_a.
\]
In particular, notice that $I$ contains the am-closed ideal $_aI$ and is contained in the am-open ideal $I_a$. Several
problems we encountered early in the study of am-ideals required us to find an answer to the dual questions: 
Is there a largest am-closed ideal contained in $I$, and is there a smallest am-open ideal containing $I$? 
Both questions have an affirmative answer which sheds interesting light on the structure of am-ideals.

\subsection{\leftline{Am-closed ideals}}\label{s:closed}
For the first question, a moment's reflection shows that it is equivalent to asking whether the sum of two am-closed ideals is am-closed. 
To answer it, we needed some tools from majorization theory (see also Section \ref {S: Majorization}). 
From the definition of $I^-$ we see that
 \begin{align*}
\xi\in \Sigma(I^-)&=\Sigma(_a(I_a)) \Leftrightarrow \xi_a \in \Sigma(I_a)
\Leftrightarrow \xi_a \le \eta_a \text{ for some }\eta\in \Sigma(I)\\
&\Leftrightarrow \text{ for every }n, \sum_{j=1}^n\xi_j \le \sum_{j=1}^n\eta_j  \text{ for some }\eta\in \Sigma(I)
\end{align*}
The last condition is called majorization and we denote it by $\xi \prec \eta$ 
(for finite sequences it is more common  to call it weak majorization and to use a different notation). 
Thus an ideal is am-closed if and only if it is hereditary (i.e., solid) under the majorization order.  
Markus \cite [Lemma 3.1]{aM64} proved that $\xi \prec \eta$ if and only if there is a substochastic matrix $P$ (i.e., a matrix with non-negative entries with row and column sums bounded by $1$) for which $\xi=P\eta$. This was the key tool in proving the theorem:
\bT{T:sum closed}\cite [Theorem 2.9]{vKgW04-Soft}  $(I+J)^- = I^- + J^-$ for all ideals $I,J$. In particular, the sum of two am-closed ideals is am-closed.
\eT

As a consequence, the collection of all the am-closed ideals contained in an ideal $I$ is directed and hence its union is an ideal which we symbolize by
\[
I_- := \bigcup\{J \mid J \subset I\text{ and }J \text{ is
an am-closed ideal}\}.
\]
$I_-$ can be easily shown to be am-closed, and hence it is the largest am-closed ideal contained in $I$. 
When $I$ is countably generated we can say more:

\bT{T:coun gen} \cite[Theorem 2.9]{vKgW04-Soft} If $I$ is a countably generated ideal, then $I_- =\, _aI$.
\eT

This theorem provides a quick and  elementary proof of a result that required some effort even in the case of principal ideals \cite[Theorem 3.11]{AS78}:

\bT {T:coun gen closed} \cite [Theorem 2.11]{vKgW04-Soft} A countably generated ideal is am-closed if and only if it is am-stable.
\eT

\subsection{\leftline{Am-open ideals}}\label{s:open}\quad
We were led to a completely different path while exploring the question of whether there is a smallest am-open ideal $I^{oo}$ containing a given ideal $I$ or, equivalently, whether the intersection of a family of am-open ideals is necessarily am-open. 
These questions arise naturally in trying to identify the am-interior of an ideal. 
From the definition of $I^o$ we see that $\xi\in \Sigma(I^o)=\Sigma((_aI)_a)$ if and only if $\xi\le \eta_a$ for some $\eta_a\in  \Sigma(I)$. 
This led us to ask whether there is a smallest average that dominates $\xi$ and hence a smallest am-open ideal that contains $(\xi)$. 
The answer is yes and to provide it we need to employ the upper and lower monotone nondecreasing and monotone nonincreasing envelopes of a real-valued sequence $\phi$ (see \cite[Section 2.3]{DFWW}):
\[
\und \phi :=\! \Big\langle\max_{ i\leq n} \phi_i\Big\rangle, \quad 
\lnd \phi :=\! \Big\langle \inf_{i\geq n} \phi_i\Big\rangle, \quad 
\uni  \phi :=\! \Big\langle \sup_{i\geq n} \phi_i\Big\rangle, \quad
\lni  \phi :=\! \Big\langle \min_{i\leq n} \phi_i\Big\rangle.
\]
\begin{lemma} \label{lem:2.14}\cite[Lemma 2.14]{vKgW04-Soft} For every $\xi  \in  c_o^*$,
\begin{enumerate}
\item[(i)]
$(\xi )^o = (\omega \,   \lnd  \frac{\xi}{\omega})$

\item[(ii)]
$(\omega \und  \frac{\xi}{\omega})$ is the smallest am-open ideal containing $(\xi )$. 
\end{enumerate}
\end{lemma}

In fact, $\omega\und \frac{\xi}{\omega} \asymp \omega \psi$ where $\psi$  is the smallest concave sequence that dominates $\und \frac{\xi}{\omega}$ and then it follows that $\omega \psi$ is the smallest  average dominating $\xi$. In contrast, it is true that $\omega\lnd \frac{\xi}{\omega}\asymp \omega \phi$ where $\phi$  is the smallest concave sequence that dominates $\frac{1}{2} \lnd \frac{\xi}{\omega}$ and that $\omega \phi$ is an average and $\omega \phi  \le \xi$, but  unless $\xi$ is itself an average, $\omega \phi$ is smallest only in the O-sense.

This lemma leads to a number of useful consequences.
\bT{T:open}  \cite[Section 2] {vKgW04-Soft} For every ideal $I$:
\item[(i)] $\Sigma(I^o) = \{\xi  \in  c_o^* \mid \omega \und \frac{\xi}{\omega} \in
\Sigma(I)\} = \{\xi  \in  c_o^* \mid \xi  \leq \omega  \lnd
\frac{\eta}{\omega}$ for some $\eta \in  \Sigma(I)\}$.
\item[(ii)] Intersections of am-open ideals are am-open.
\item[(iii)] $I^{oo}:=
\bigcap\{J \mid J \supset I \text{ and $J$ is am-open} \}$
is the smallest am-open ideal containing $I$.
\item[(iv)] $
\Sigma(I^{oo}) = \{\xi  \in  c_o^* \mid \xi  \leq \omega \und \frac{\eta}{\omega} \text{ for some  } \eta \in  \Sigma(I)\}$.
\item[(v)] For all $I,J$, $(I + J)^{oo} = I^{oo} + J^{oo}$ and  $I^o + J^o \subset (I + J)^o$.
The inclusion can be proper.
\eT

\subsection{\leftline{\bf Arithmetic mean at infinity ideals}}\label{S: am infty ideals}\quad
For an  ideal $I$ contained in the trace-class, the first order arithmetic mean ideals are always trivial: $_aI=I^o={0}$, $I^-=\mathscr L_1$, and $I_a=(\omega)$ but there is a theory of  am-$\infty$ ideals that somewhat mirrors the theory of am-ideals. 

\bT{T:aminfty closed} \cite[Section 3]{vKgW04-Soft}
\item[(i)] $(I + J)^{-\infty}  = I^{-\infty} + J^{-\infty}$ for all ideals $I$, $J$. 
\item[(ii)] The sum of two am-$\infty$ closed ideals is am-$\infty$ closed.
\item[(iii)] $I_{-\infty}  := \bigcup\,\{J \mid J \subset I$  and $J$ is am-$\infty$ closed is the largest am-closed ideal contained in $I$.
\item[(iv)] If $I$ is a countably generated ideal, then $I_{-\infty} = \text{$_{a_\infty}$}I$.
\item[(v)]  A countably generated ideal is am-$\infty$ closed if and only if it is am-$\infty$ stable.
\end{theorem}
 
 While the statements are similar to those of the am-case, the techniques employed in the proofs are substantially different. For instance, the proof of (i) depends on a $w^*$-compactness argument rather than the Markus based matricial one used in Theorem 1.5.\\
 
 The characterization of am-$\infty$ open ideals depends on the following lemma:
\begin{lemma}\label{lem:3.7} \cite[Lemma 3.7]{vKgW04-Soft}
For every principal ideal $I$, the following are equivalent.
\item[(i)] $I$ is am-$\infty$ open.
\item[(ii)] $I = (\eta _{a_\infty})$ for some $\eta \in (\ell^1)^*$.
\item[(iii)] $I = (\xi )$ for some $\xi$ for which $\frac{\xi}{\omega}  \in c_o^*$.
\end{lemma}

 The crux of  the proof  of (iii) $\Rightarrow$ (ii) consist in showing that if $\psi$ is the largest convex sequence pointwise dominated by $\frac{\xi}{\omega}$, then $\frac{1}{2}D_{\frac{1}{2}}\frac{\xi}{\omega}\le \psi$. 

\bT{T:aminfty open} \cite[Section 3]{vKgW04-Soft}
For every ideal $I$:
\item[(i)] $\Sigma(I^{o\infty}) = \{\xi  \in  \Sigma( \se (\omega))
\mid \omega  \uni \frac{\xi}{\omega} \in \Sigma(I)\} 
= \{\xi  \in  c_o^* \mid \xi \leq \omega \lni \frac{\eta}{\omega} \text{ for some } 
\eta \in \Sigma(I \cap \se (\omega))\}$ 
where the operation $\se (\omega)$ is defined below.
\item[(ii)] Intersections of  am-$\infty$ open ideals are am-$\infty$ open.
\item[(iii)] $I^{oo\infty} := \bigcap\{J \mid   I \cap \se (\omega) \subset J \}$ is the smallest am-$\infty$ open ideal containing $I$.
\item[(iv)] $\Sigma(I^{oo\infty}) = \big\{\xi  \in c_o^* \mid \eta \leq \omega
\uni \frac{\eta}{\omega} \text{ for some }\eta \in \Sigma(I \cap \se (\omega))\big\}$.
\item[(v)] $I^{oo\infty}+ J^{oo\infty}  = (I +J)^{oo\infty}$ and $I^{o\infty}+ J^{o\infty}  \subset (I +
J)^{o\infty}$ for any two ideals $I$ and $J$ and the inclusion can be proper.
\eT

 A further substantial difference from the am-case is due to the fact that, while for any sequence $\xi\in \co*$ the sequence $\xi_a$ always satisfies the $\Delta_{1/2}$ condition and also  $\xi \le \xi_a$, in the am-$\infty$ case neither of these facts hold (see \cite [Example 4.5]{vKgW04-Traces}). 
In fact the two conditions are related.
 
  \bP{P:am-infty seq}\cite [Corollary 4.4] {vKgW04-Traces}
Let $\xi \in (\ell^1)^*$. 
\item[(i)] If $\xi$ satisfies the $\Delta_{1/2}$-condition, so does $\xi_{a_\infty}$.
\item[(ii)] $\xi_{a_\infty}$ satisfies the $\Delta_{1/2}$-condition if and only if $\xi = O(\xi_{a_\infty})$.
\eP

It is not surprising and relatively simple to show that 
if $\xi\in (\ell^1)^*$, then $(\xi)\subset (\xi)_{a_\infty} =  (\xi_{a_\infty})$ \cite[Lemma 4.7]{vKgW04-Traces}, 
and hence if $\xi_{a_\infty}= O (\xi)$ then $(\xi) = (\xi)_{a_\infty}$, i.e., $(\xi) $ is am-$\infty$ stable. 
The opposite implication, however, is non-trivial. 
By analogy with the am-case, call regular at infinity a sequence $\xi \in (\ell^1)^*$ for which $(\xi)=(\xi)_{a_\infty}$.

\bT {T:infty regular}\cite [Theorem 4.12] {vKgW04-Traces}
If $\xi \in (\ell^1)^*$, then the following conditions are equivalent.
\item[(i)]   $\xi$ is regular at infinity.
\item[(ii)]  $\xi_{a_\infty} = O(\xi)$
\item[(iii)] There are constants $C > 0$ and $p > 1$ for which $\xi_n \leq C (\frac{m}{n})^p \xi_m$ for all integers $n \geq m$.
\item[(iv)]  $\xi_{a_\infty}$ is regular at infinity.
\eT
  
\noindent Explicitly, the principal ideal $(\xi)$ is am-$\infty$ stable if and only if $\xi_{a_\infty} = O(\xi)$.

  Condition (ii) strengthens \cite[Corollary 5.19]{DFWW} and \cite[Corollary 5.6]{mW02} by eliminating the need for ampliations, i.e., replacing the condition $\xi_{a_\infty} = O(D_m\xi)$ for some $m$ by the condition $\xi_{a_\infty} = O(\xi)$.
  
Condition (iii) says that $\alpha (\xi)< -1$ where $\alpha (\xi)$ is the Matuszeswka index of the sequence $\xi$  (see \cite{DFWW}) which can be defined as
 \[ 
  \alpha(\xi) = \inf \{\gamma \mid \exists\,C>0~\text{such that}~ \xi_n \leq C\left(\frac{n}{m}\right)^\gamma \xi_m ~\text{for all}~ n \geq m\}. 
 \]

\section{\leftline{\bf Softness properties of operator ideals}}\label{S:3}

\subsection{\leftline{\bf Definitions}}\label{s: Definitions} \quad
In \cite[Section 4]{vKgW04-Soft} we studied \textit{soft ideals} first introduced in \cite[Section 4]{vKgW02} and \cite[Section 4]{vKgW04-Traces}, motivated in part by Dixmier's implicit use of softness to construct the so-called Dixmier trace \cite{jD66}. 

\begin{definition}\label{def:4.1}
The soft interior of an ideal $I$ is the product $\se I: = IK(H)$,
i.e., the ideal with characteristic set
\[
\Sigma(\se I) = \{\xi  \in  c_o^* \mid \xi  \leq \alpha \eta
\text{ for some } \alpha \in  c_o^*,\ \eta \in  \Sigma(I)\}.
\]
The soft cover of an ideal $I$ is the quotient $\scop I:= I : K(H)$,
i.e., the ideal with characteristic set
\[
\Sigma(\scop I) = \{\xi  \in  c_o^* \mid \alpha\xi  \in  \Sigma(I)
\text{ for all } \alpha \in  c_o^*\}.
\]
An ideal is called soft-edged if $\se I = I$ and soft-complemented
if $\scop I = I$. 

\noindent A pair of ideals $I \subset J$ is called a soft pair if $\se J = I$
and $\scop I = J$.
\end{definition}

\noindent Soft-complemented and soft-edged ideals share many common properties. Collectively, we call them \textit{soft ideals}. The terminology we introduced is motivated by the fact that $I$ is soft-edged if and only if, for every $\xi \in \Sigma(I)$, 
one has $\xi = o(\eta)$ for some $\eta \in  \Sigma(I)$. 
Analogously, an ideal $I$ is soft-complemented if and only if, for every $\xi  \in  c_o^* \setminus \Sigma(I)$, 
one has $\eta = o(\xi)$ for some $\eta \in c_o^* \setminus \Sigma(I)$.

The soft-interior $\se I$ and the soft-cover $\scop I$ are, respectively, 
the largest soft-edged ideal contained in $I$ and the smallest soft-complemented ideal containing $I$. 
The pair $\se I \subset \scop I$ is the generic example of what we call a soft pair. 
Many classical ideals, i.e., ideals whose characteristic set generates a classical sequence space, turn out to be soft and form natural soft pairs 
as discussed below. 

To prove soft-complementedness of an ideal we often found it convenient to prove instead a stronger property which we call strong soft-complementedness 
\cite[Proposition 4.5]{vKgW04-Soft}. 

\bD{D:4.4}
An ideal $I$ is said to be strongly soft-complemented (ssc for short) if
for every countable collection of sequences $\{\eta^{(k)}\} \subset
c_o^* \setminus \Sigma(I)$ there is a sequence of indices $n_k \in
\mathbb{N}$ such that $\xi  \not\in \Sigma(I)$ whenever $\xi  \in  c_o^*$ and
$\xi_i \geq \eta_i^{(k)}$ for all $k$ and for all $1 \leq i \leq n_k$.
\eD

\vspace{0.1cm}

\subsection{\leftline{\bf Classical ideals}}\label{s: Classical ideals}

\noindent\textbf{Countably generated ideals.} 
\bP{P: Countably generated ideals-softness}\cite [Proposition 4.6]{vKgW04-Soft} 
\item[(i)] Countably generated ideals, in particular, principal ideals, are strongly soft-complemented and hence soft-complemented.
\item[(ii)] A principal ideal $(\xi)$ is soft-edged if and only if $\xi = o(D_m \xi)$ for some $m$.
\eP

\noindent\textbf{Banach ideals} are ideals that are complete with respect to a symmetric norm (see for instance \cite[Section 4.5]{DFWW})
and were called uniform-cross-norm ideals by Schatten \cite{rS60},
symmetrically normed ideals by  Gohberg and Krein  \cite{GK69}, and
symmetric norm ideals by other authors. Recall that the norm of
$I$ induces on the finite rank ideal $F$ (or, more precisely, on $S(F)$,
the associated space of sequences of $c_o$ with finite support) a
symmetric norming function $\phi$, and the latter permits one to
construct the so-called minimal Banach ideal
$\mathfrak{S}_\phi^{(o)} = \cl(F)$ contained in $I$ (the closure taken
in the norm of $I$) and  the maximal Banach ideal $\mathfrak{S}_\phi$ containing $I$ where
\begin{align*}
\Sigma\big(\mathfrak{S}_\phi\big) &= \big\{\xi  \in  c_o^* \mid \phi(\xi):=\, \sup \phi\big( \langle \xi_1, \xi_2, \dots , \xi_n, 0, 0, \ldots \rangle \big) < \infty \big\} \\
\Sigma\big(\mathfrak{S}_\phi^{(o)}\big) &= \big\{\xi \in \Sigma\big(\mathfrak{S}_\phi\big) \mid \phi\big(\langle \xi_n, \xi_{n+1}, \ldots \rangle \big) \rightarrow 0\big\}.
\end{align*}
If $I$ is a Banach ideal, the ideals $\mathfrak{S}_\phi^{(o)}$ and $\mathfrak{S}_\phi$ can be obtained from $I$ through a ``soft'' operation, i.e., 

\bP{P:Gohberg/Krein soft pairs}\cite[Proposition 4.7]{vKgW04-Soft} \\
If $I$ is a Banach ideal, then
$\mathfrak{S}_\phi^{(o)} = \se I$ and $\mathfrak{S}_\phi = \scop I$, 
the embedding $\mathfrak{S}_\phi^{(o)} \subset \mathfrak{S}_\phi$ is a natural example of a soft pair, and $\mathfrak{S}_\phi$ is ssc. In particular, $\se I$ and $\scop I$ are also Banach ideals.
\eP

Recall that in the notations of \cite{DFWW} and \cite{vKgW04-Soft}, Gohberg and Krein \cite{GK69} showed that the symmetric norming function $\phi(\eta) := \sup \frac{\eta_a}{\xi_a}$ induces a complete norm on the am-closure $(\xi )^-$ of the principal ideal $(\xi)$ and for this norm,
\[
\cl(F) = \mathfrak{S}_\phi^{(o)} \subset  \cl (\xi ) \subset
\mathfrak{S}_\phi = (\xi)^-.
\]
The fact that $\mathfrak{S}_\phi$ is soft-complemented was obtained in \cite[Theorem 3.8]{nS74}, 
but Salinas proved only that (in our notations) $\se \mathfrak{S}_\phi \subset \mathfrak{S}_\phi^{(o)}$ \cite[Remark 3.9]{nS74}. 
Varga reached the same conclusion in the case of the am-closure of a principal ideal with a non-trace class generator \cite[Remark 3]{jV89}.\\

The fact that $\mathfrak{S}_\phi^{(o)} \subset \mathfrak{S}_\phi$ is always a soft pair yields immediately the equivalence of parts (a)--(c) in \cite[Theorem 2.3]{nS74} without considering norms and hence establishes (d) and (e).\\
The same fact  can also be used to obtain an elementary proof that $(\xi)=\cl(\xi)$ if and only if $\xi$ is regular \cite[Proposition 4.9]{vKgW04-Soft}. This was proven by Allen and Shen \cite[Theorem 3.23]{AS78} using Salinas' results \cite{nS74} on (second) K\"{o}the duals and by Varga \cite[Theorem 3]{jV89} and is a special case of \cite[Theorem 2.36]{DFWW}.

The fact that $\mathfrak{S}_\phi^{(o)}= \se \mathfrak{S}_\phi$ brings a further  consequence.  As remarked in  \cite[Section 4.9]{DFWW},  $ \mathfrak{S}_\phi$ is am-closed. By Proposition \ref {P;other comm} below, $\mathfrak{S}_\phi^{(o)}$ is also am-closed.\\

\noindent\textbf{Orlicz ideals.} 
Recall from \cite[Sections 2.37 and 4.7]{DFWW} that if $M$ is a monotone nondecreasing function on $[0,\infty)$ with $M(0) = 0$, then the small Orlicz ideal
$\mathscr L_M^{(o)}$ is the ideal with characteristic set $\{\xi \in  c_o^* \mid \sum_n M(t\xi_n) < \infty \text{ for all } t > 0\}$ and
the Orlicz ideal $\mathscr L_M$ is the ideal with characteristic set $\{\xi  \in  c_o^* \mid \sum_n M(t\xi_n) < \infty \text{ for some } t > 0\}$. 

\begin{proposition}\label{prop:4.11}\cite[Proposition 4.11]{vKgW04-Soft}
Let $M$ be a monotone nondecreasing function on $[0, \infty)$ with $M(0) = 0$. 
Then $\mathscr L_M^{(o)}$ is soft-edged, $\mathscr L_M$ is ssc, and $\mathscr L_M^{(o)} \subset \mathscr L_M$ is a soft pair.
\end{proposition}

\noindent The fact that $\mathscr L_M^{(o)} \subset  \mathscr L_M$ forms a soft pair can simplify proofs of some properties of Orlicz ideals.
Indeed, together with  \cite[Proposition 3.4]{vKgW04-Traces} that states that for an ideal $I$, $\se I$ is am-stable if and only if $\scop I$ is am-stable 
if and only if $I_a \subset  \scop I$, and combined with \cite[Lemma 4.16]{vKgW04-Soft} it yields an immediate proof of the following results in \cite{DFWW}: the equivalence of (a), (b), (c) in Theorem 4.21 and hence the equivalence of (a), (b), (c) in Theorem 6.25, the equivalence of (b), (c), and (d) in
Corollary 2.39, the equivalence of (b) and (c) in Corollary 2.40, and the equivalence of (a), (b), and (c) in Theorem 3.21.\\

\noindent\textbf{Lorentz ideals.} If $\phi$ is a monotone nondecreasing nonnegative sequence satisfying the $\Delta_2$-condition, 
i.e., $\sup\frac{\phi_{2n}}{\phi_n} < \infty$, then in the notations of \cite[Sections 2.25 and 4.7]{DFWW} the Lorentz ideal $\mathscr L(\phi)$ corresponding to the sequence space $\ell(\phi)$ is the ideal with characteristic set
\[
\Sigma(\mathscr L(\phi)) := \,\bigg\{\xi  \in  c_o^* \mid \|\xi \|_{\ell(\phi)} := \sum_n \xi_n(\phi_{n+1} - \phi_n) < \infty\bigg\}.
\]
A special case of Lorentz ideal is the trace class $\mathscr L_1$ which corresponds to the sequence $\phi = \langle n \rangle$ and the sequence space 
$\ell(\phi) = \ell^1$.

\begin{proposition}\label{prop:4.12}\cite[Proposition 4.12]{vKgW04-Soft}
If $\phi$ be a monotone nondecreasing nonnegative sequence satisfying the $\Delta_2$-condition, then $\mathscr L(\phi)$ is both soft-edged and strongly
soft-complemented.
\end{proposition}
\noindent In particular, $\mathscr L_1$ is both soft-edged and soft-complemented.\\

\noindent\textbf{K\"{o}the duals}
\begin{proposition}\label{prop:4.13}\cite[Proposition 4.13]{vKgW04-Soft}
Let $I$ be a soft-complemented ideal and let $X$ be a nonempty subset of $[0,\!\infty)^{\mathbb N}$. 
Then the ideal with characteristic set $\Sigma(I) : X$ is soft-complemented.
In particular, K\"{o}the duals are soft-complemented.
\end{proposition}

The K\"{o}the dual, $I := \{\langle e^n \rangle\}^\times$ of the singleton $\{\langle e^n \rangle\}$, is soft-complemented but is not strongly soft-complemented \cite[Example 4.15]{vKgW04-Soft}.\\

\noindent\textbf{Idempotent ideals} are ideals for which $I = I^2$. 
Notice that an ideal is idempotent if and only if $I = I^p$ for some positive $p \neq 0,1$, if and only if $I = I^p$ for all $p > 0$.

\begin{proposition}\label{prop:4.17}\cite[Proposition 4.17]{vKgW04-Soft}
Idempotent ideals are soft-edged and soft-complemented.
\end{proposition}

\noindent The remarks following \cite[Proposition 5.3]{vKgW04-Soft} show that idempotent ideals may fail to be strongly soft-complemented.\\

\noindent\textbf{Marcinkiewicz ideals} are the pre-arithmetic means of principal ideals, and we consider also their am-$\infty$ analogs. 
That these ideals are strongly soft-complemented follows from the following proposition combined with 
Proposition \ref{P: Countably generated ideals-softness}(i).

\begin{proposition} \label{prop:4.18}\cite[Proposition 4.18]{vKgW04-Soft}
The pre-arithmetic mean and the pre-arithmetic mean at infinity of a strongly soft-complemented ideal are strongly soft-complemented. 

In particular, Marcinkiewicz ideals are strongly soft-complemented.
\end{proposition}

\subsection{\leftline{\bf Operations on soft ideals}}\label{s: soft}\quad
In applying ``softness" properties of operator ideals it is often essential to know which operations preserve softness. 
In some cases, the answer comes from an elementary computation, e.g.,

\bP{P:inters}\cite [Proposition 5.1, Example 5.2] {vKgW04-Soft}\\ 
For every collection of ideals $\{I_\gamma , \gamma \in \Gamma\}$:
\item[(i)]
 $\bigcap_\gamma  \se I_\gamma \supset \se (\bigcap_\gamma  I_\gamma)$. If $\Gamma$ is infinite, the inclusion can be proper.
\item[(ii)]
 $\bigcap_\gamma  \scop I_\gamma = \scop (\bigcap_\gamma  I_\gamma)$ \\
In particular, the intersection of soft-complemented ideals is soft-complemented. However, the intersection of an infinite countable strictly decreasing chain
of principal ideals is never strongly soft-complemented \cite [Proposition 5.3]{vKgW04-Soft}
\end{proposition}

Even more elementary is the fact that the $\se$ operation distributes for sums: $\se(I+J)= K(H)(I+J)= \se I + \se I$. 
Thus we found surprisingly challenging the analogous question about the ``sc" operation. We were only able to obtain:

\bT{T:sc sum}\cite [Theorem 5.7] {vKgW04-Soft} The sum of an ssc ideal  and a countably generated ideal
 is ssc and hence soft-complemented.
\eT

\noindent We do not know if these conditions can be relaxed:\\ 

\noindent \textbf{Question 1.}
Is the sum of two ssc ideals ssc? Is it soft-complemented? Is the sum of a soft-complemented ideal  and a countably generated ideal soft-complemented?\\

In dealing with the codimension question for traces (see next section) we heavily depended on commutation relations between ``soft operations" and arithmetic mean operations and  arithmetic mean at infinity operations. 

\bT{T:comm}\cite [Theorem 6.1]{vKgW04-Soft} 
Let $I$ be an ideal.

\item[(i)] $\scop \text{$_a$}I \subset \text{$_a$}(\scop I)$

\item[(i$'$)] $\scop \text{$_a$}I = \text{$_a$}(\scop I) \text{ if and only if }
\omega \not\in \Sigma(\scop I) \setminus \Sigma(I)$

\item[(ii)] $\se I_a \subset (\se I)_a$

\item[(ii$'$)] $\se I_a = (\se I)_a\text{ if and only if }I = \{0\} \text{ or }I
\not\subset \mathcal{L}_1$

\item[(iii)] $\scop I_a \supset (\scop I)_a$

\item[(iv)] $\se \text{$_a$}I \supset \text{$_a$}(\se I)$

\item[(iv$'$)] $\se \text{$_a$}I = \text{$_a$}(\se I) \text{ if and only if }
\omega \not\in \Sigma(I) \setminus \Sigma(\se I)$.
\eT
 As a consequence, and essential for the next section is the fact that  $\se I$ is am-stable if and only if so is $\scop I$ and both are implied by (but not equivalent to) the am-stability of $I$ and both can be  linked to the following condition involving the commutator space of $I$:

 \bP{P: se stability} \cite [Proposition 4.20, 3.4]{vKgW04-Traces}  Let $I$ be an ideal. Then $\se I \subset F + [I,B(H)]$ if and only if
 \item[(i)] $\se I$ is am-$\infty$ stable in the case that $\omega \notin \Sigma(I)$
 \item[(ii)] $\se I$ is am-stable in the case that $\omega \in \Sigma(I)$\\
Also, $\se I$ is am-$\infty$ stable (resp., am-stable) if and only if $\sc I$ is am-$\infty$ stable (resp., am-stable).
 \eP

The ``missing'' reverse inclusion of Theorem \ref {T:comm}(iii) is an open question and it is equivalent to a number of other questions listed below:\\ 

\noindent \textbf{Question 2.}  Let $I$ be an ideal.
\begin{enumerate}
\item[(i)]  Is (or under what condition(s)) $\scop I_a = (\scop I)_a$?
\item[(i$'$)]  Is (or under what condition(s)) $\scop I^- = (\scop I)^-$? \\
Equivalently, is the am-closure of a soft-complemented ideal also soft-complemented?
 \item[(ii)] Is (or under what condition(s)) $ (\scop I)_a$ soft-complemented?
  \item[(ii$'$)] Is (or under what condition(s)) $ (\scop I)^-$ soft-complemented?\\
\end{enumerate}

  Further results  linking softness properties and arithmetic mean properties are
  
  \bP{P;other comm} \cite [Corollaries 6.7, 6.5] {vKgW04-Soft} 
  \item[(i)] Let $I$ be soft-complemented. Then so are $~_aI$, $I^o$, and
$I_-$.  Let $I$ be soft-edged. Then so are $~_aI$, $I^o$, and $I^-$; $I_a$ is soft-edged if and only if
either $I = \{0\}$ or $I \not\subset \mathcal{L}_1$; $I^{oo}$ is soft-edged if and only
if either $I = \{0\}$ or $I \not\subset (\omega)$.
 \item[(ii)] Let $I$ be am-closed. Then so are $\scop I$ and $\se I$. Let $I$ be am-open. Then $\scop I$ is am-open while $\se I$ is am-open if and only if $I\neq(\omega)$.
 \eP

Similar results and similar open question (but as usual with somewhat more complicated conditions and proofs) hold for arithmetic mean at infinity operations (see \cite [Section 6]{vKgW04-Soft}).

\section{\leftline{\bf Traces and the codimension conjecture }}\label{S:4}

The first major question we  investigated in our study of arithmetic mean operator ideals in \cite{vKgW04-Traces} was :
``How many nonzero traces can an ideal support?" Because traces on an ideal $I$ are the unitarily invariant linear functionals on $I$ and their unitary invariance is equivalent to their vanishing on $[I,B(H)]$, they can be identified with the elements of the linear dual of the quotient space $\frac{I}{[I,B(H)]}$. Our work in \cite{vKgW04-Traces} led us to the following question:\\

\noindent \textbf{Question 3.}
Is the following true for all ideals $I$?
\[
\di \frac{I}{[I,B(H)]} ~\in
\begin{cases}
\{1,\infty\}    &\text{when}~\omega \notin \Sigma(I) \\
\{0,\infty\}    &\text{when}~\omega \in \Sigma(I).
\end{cases}
\]\\

We conjecture that the answer to this question is affirmative and we will present in Section \ref {S:2.2.2} some evidence to support this conjecture. 
It was precisely to  investigate this conjecture and related questions that we developed the tools presented in the two preceding sections.

\subsection{\leftline{\bf Trace extensions, hereditariness, and the uniqueness of traces}}\label{S:4.1}
It is well-known that
the restriction of a trace on a nonzero ideal $I$ to the ideal $F$ of finite rank operators
must be a (possibly zero) scalar multiple of the standard trace $Tr$.
In other words, every trace on a nonzero ideal is a trace extension of some scalar multiple of $Tr$ on $F$.

\begin{definition}\label{D:singular trace} A trace vanishing on $F$ is called singular, and nonsingular otherwise.
\end{definition}

Dixmier \cite{jD66} provided the first example of a (positive) singular trace.
Its domain is the
am-closure, $(\eta)^-=\,_a(\eta_a)$, of a principal ideal
$(\eta) \subset \se(\eta)_a$.

Theorem \ref{T:DFWW} yields a complete characterization of ideals that support a nonsingular trace,
namely, those ideals that do not contain $\diag \omega$
(cf. \cite[Introduction, Application 3 of Theorem 5.6]{DFWW} and also \cite{kDgWmW00}).
Generalizing this argument proves the following:

\begin{proposition}\label{P:trace extension} \cite[Proposition 5.3]{vKgW04-Traces} 
Let $I$ and $J$ be ideals and let $\tau$ be a trace on $J$.
Then $\tau$ has a trace extension to $I+J$ if and only if
\[
 J\cap[I,B(H)]\subset \{\,X \in J\mid\tau(X)=0\,\}.
\]
Moreover, the extension is unique if and only if $I \subset J + [I,B(H)]$.

In particular, $\omega\notin\Sigma(I)$
if and only if $Tr$ extends from $\mathscr L_1$ to $\mathscr L_1+\,I$ if and only if $Tr$ extends from $F$ to $I$.
\end{proposition}

Further consequences of  \cite{DFWW} (see Theorem \ref{T:pos part + infty characterization} above) are

\begin{proposition}\label{P:(L1+[])+=L1+}\cite[Proposition 5.5]{vKgW04-Traces} 
Let $I$ be an ideal for which $\omega \notin \Sigma(\se I)$. Then
\[
(\mathscr L_1 + [I,B(H)])^+ = \mathscr L_1^+.
\]
In particular, $(\mathscr L_1 + [I,B(H)])^+$ is hereditary for every ideal $I$.
\end{proposition}

\begin{corollary} \label{C:(F+[])+}\cite[Corollary 6.2]{vKgW04-Traces}
 Let $I\ne \{0\}$ be an ideal.
\item[(i)] If $\omega \notin \Sigma(I)$, then $(F + [I,B(H)])^+ = (_{a_\infty}I)^+$.
\item[(ii)] If $\omega \in \Sigma(I)$, then $(F + [I,B(H)])^+ = (_aI)^+$.
\item[(iii)] $(F + [I,B(H)])^+$ is hereditary (i.e., solid).
\end{corollary}

An ideal $I$ supports a unique nonzero trace (up to scalar multiplication)
precisely when $\di \frac{I}{[I,B(H)]} = 1$. All the ideals for which $\omega \notin \Sigma(I)$
support some non-singular trace. The next theorem which is a consequence of Corollary \ref {C:(F+[])+}  tells us that this trace is unique precisely when the ideal is am-$\infty$ stable.

\begin{theorem}\label{T:i-iii}\cite [Theorem 6.6]{vKgW04-Traces}
If $I\ne \{0\}$ is an ideal where $\omega \notin \Sigma(I)$, then the following are equivalent.
\item[(i)]  $I$ supports a nonzero trace unique up to scalar multiples.
\item[(ii)]  $I \subset \mathscr L_1$ and every trace on $I$ is a scalar multiple of $Tr$.
\item[(iii)] $\di \frac{I}{[I,B(H)]} = 1$
\item[(iv)] $I = F + [I,B(H)]$
\item[(v)]  $I \subset \mathscr L_1$ and $[I,B(H)]=\{X\in I \mid Tr\,X=0\}$
\item[(vi)]  $I$ is am-$\infty$ stable, i.e., $I=~_{a_\infty}I$.
\end{theorem}

What if  $\omega \in \Sigma(I)$? We know that if $I$ is am-stable, then it supports no nonzero trace. But could it support a unique trace? Equivalently:\\

\noindent \textbf{Question 4.}
If $(\omega) \subset I \ne [I,B(H)]$, must $\di \frac{I}{[I,B(H)]} > 1$?\\

As mentioned in the introduction, a principal ideal $(\xi)$ supports \textit{no nonzero trace} precisely when $\xi$ is \textit{regular} and Theorem \ref {T:i-iii} shows that it supports a \textit{unique nonzero trace} precisely when $(\xi)= ~_{a_\infty}(\xi)$, or equivalently, when $(\xi)=(\xi) _{a_\infty}= (\xi) _{a_\infty}$, which  by Theorem \ref {T:coun gen closed}, is equivalent to $\xi_{a_\infty}= O(\xi)$. A similar characterization  of principal ideals supporting a unique trace was obtained in \cite[Corollary 5.6]{mW02}.

So far, we know that $[I,B(H)]^+$, $(F+[I,B(H)])^+$, and $(\mathscr L_1+[I,B(H)])^+$ are always hereditary. Of course, if $J$ is an ideal that either contains or is contained in the commutator space
$[I,B(H)]$, then $(J+[I,B(H)])^+$ must be also hereditary.\\

\noindent \textbf{Question 5.} Is $(J+[I,B(H)])^+$ hereditary for all pairs of ideals $I$ and $J$?
If not, when is it?\\

\subsection{\leftline{\bf Applications to elementary operators, and problems of Shulman}}\label{S:4.2}Propositions \ref {P:trace extension}  and \ref {P:(L1+[])+=L1+} find natural applications to questions on \textit{elementary operators}. \\
If $A_i, B_i \in B(H)$, then the $B(H)$-map
\[
B(H) \owns T \rightarrow \Delta(T) := \sum_{i=1}^{n} A_iTB_i
\]
is called an elementary operator and its adjoint $B(H)$-map is $\Delta^*(T) := \sum_{i=1}^{n} A_i^*TB_i^*$. \linebreak
Elementary operators include commutators and intertwiners and hence their theory is connected to the structure of commutator spaces.
The Fuglede-Putnam Theorem \cite{bF51}, \cite{cP51} states that for the case $\Delta(T) = AT-TB$ where A, B are normal operators,
$\Delta(T) = 0$ implies that $\Delta^*(T) = 0$.
Also for $n = 2$, Weiss \cite{gW83}
generalized this further to the case where \{$A_i$\} and \{$B_i$\}, $i=1,2$, are
separately commuting families of normal operators by proving that
$\Delta(T) \in \mathscr{L}_2$ implies $\Delta^*(T) \in \mathscr{L}_2$ and
$\Vert\Delta(T)\Vert_2 = \Vert\Delta^*(T)\Vert_2$.
(This is also a consequence of Voiculescu's
\cite[Theorem 4.2 and Introduction to Section 4]{dV79/81} but neither Weiss' nor Voiculescu's methods seem to apply to the case $n > 2$.)
In \cite{vS83} Shulman showed that for $n = 6$, $\Delta(T) = 0$ does not imply $\Delta^*(T) \in \mathscr{L}_2$.

If we impose some additional conditions  on the families $\{A_i\},\{B_i\}$ and $T$,
we can extend these implications to arbitrary $n$ past the obstruction found by Shulman.

Assume that \{$A_i$\}, \{$B_i$\}, $i = 1,\dots,n$,
are separately commuting families of normal operators and let $T \in B(H)$.

Define the following ideals:
\begin{align*}
&L := (\sum_{i=1}^{n}(A_i T)(B_i))^2, \quad L_* := (\sum_{i=1}^{n}(A_i^* T)(B_i))^2, \\
&R := (\sum_{i=1}^{n}(A_i)(TB_i))^2, \quad R_* := (\sum_{i=1}^{n}(A_i)(TB_i^*))^2 ~\text{and}\\
&I_{\Delta,T} := L \cap L_* \cap R \cap R_*, \quad S = (\sum_{i=1}^{n}(A_i T B_i)) \cap L^{1/2}\cap R^{1/2}
\end{align*}
where $(X)$ denotes the principal ideal generated by the operator $X$ and $MN$ (resp., $M+N$) denotes the product (resp., sum) of the ideals $M,N$.
Then $I_{\Delta,T}$ and $S$ are either $\{0\},~B(H)$, or a principal ideal.

\begin{proposition}\label{P:FugL2}\cite [Proposition 5.7]{vKgW04-Traces} 
If $\omega \notin \Sigma(I_{\Delta,T})$, then
$\Delta(T) \in \mathscr{L}_2$ implies \begin{center}$\Delta^*(T) \in \mathscr L_2$ and
$\Vert\Delta(T)\Vert_2 = \Vert\Delta^*(T)\Vert_2$.\end{center}
\end{proposition}

\noindent A sufficient condition independent of $T$ that insures that $\omega \notin \Sigma(I_{\Delta,T})$ is
\begin{center}$\omega^{1/4} \ne O(\sum_{i=1}^n(s(A_i)+s(B_i)))$. \end{center}
So also is the condition
$\omega^{1/2} \ne O(\sum_{i=1}^n s(A_i))$ or the condition $\omega^{1/2} \ne O(\sum_{i=1}^n s(B_i))$.\\

Another application of Propositions \ref{P:trace extension} and \ref{P:(L1+[])+=L1+} is  to a problem of Shulman. \\
Let $\Delta(T)=\sum_{i=1}^{n} A_iTB_i$ be an elementary operator where the operators $A_i$ and $B_i$ are not assumed to be  commuting or normal. Shulman showed that the composition $\Delta^*(\Delta(T)) = 0$ does not imply $\Delta(T) = 0$
and conjectured that this implication holds under the additional
assumption that $\Delta(T) \in \mathscr{L}_1$.
In the case that the ideal $S$ is ``not too large" we can prove the implication without making this assumption.

\begin{proposition}\label{P:Shulman problem}\cite[Proposition 5.8]{vKgW04-Traces} 
If $\omega \notin \Sigma(S)$, then $\Delta^*(\Delta(T)) \in \mathscr{L}_1$ implies that
$\Delta(T) \in \mathscr{L}_2$
and $\Vert\Delta(T)\Vert_2 = Tr~ T^*\Delta^*(\Delta(T))$.

In particular, if $\Delta^*(\Delta(T)) = 0$ then $\Delta(T) = 0$.
\end{proposition}

The key tools used in both propositions are the hereditariness of $(\mathscr L_1 + [I,B(H)])^+$ established in Proposition \ref {P:(L1+[])+=L1+}  and the fact that $[I,J]=[IJ,B(H)]$, which was established in \cite{DFWW}.
The latter fact was also recently similarly used in \cite{FKM06}.

\subsection{\leftline{\bf Commutator spaces with infinite codimension}}\label{S:2.2.2}\quad
In this section we present some conditions under which $[I,B(H)]$ has infinite codimension in $I$.
Our key technical result is:
\begin{theorem}\label{T: BASICS: seJ notin F+[I,B(H)] implies uncountable dimension-(J+[I,B(H)])/(F+[I,B(H)])}\cite[Theorem 7.1]{vKgW04-Traces}
If $I$ and $J$ are ideals and $\se J \not\subset F+[I,B(H)]$ then
$\frac{J+[I,B(H)]}{F+[I,B(H)]}$ has uncountable dimension.
 \end{theorem}
 So, if we take $J=I$, we see that if $\se I \not\subset F+[I,B(H)]$, then $\frac{I}{[I,B(H)]}$ has uncountable dimension. 
The condition $\se I \subset F + [I,B(H)]$ is equivalent to the am-stability of $\se I$ for the case when $\omega \in \Sigma (I)$ 
\cite [Proposition 3.4]{vKgW04-Traces} or the am-$\infty$ stability of  $\se I$ for the case when $\omega \notin \Sigma (I)$ 
\cite[Proposition 4.20]{vKgW04-Traces}.  As a consequence, recalling that  $st^a(\mathscr L_1)$ is the smallest am-stable ideal and that $st_{a_\infty}(\mathscr L_1)$ is the largest am-$\infty$ stable ideal, we obtain the following. 

\begin{theorem}\label{T: stabilizer containments and uncountable dimension}\cite[Theorem 7.2]{vKgW04-Traces}
Let $I \ne \{0\}$ be an ideal.
Then $\frac{I}{[I,B(H)]}$ has uncountable
dimension if any of the following conditions hold.
\item[(i)] $I \subset st_{a_\infty}(\mathscr L_1)$ and $\se I$ is not am-$\infty$ stable.
\item[(ii)] $I \supset st^a(\mathscr L_1)$ and $\se I$ is not am-stable.
\item[(iii)] $I \not\subset st_{a_\infty}(\mathscr L_1)$ and $I \not\supset st^a(\mathscr L_1)$.
\end{theorem}

In particular, we see that for  ``intermediate ideals" , $[I,B(H)]$ always has infinite codimension and that for  ``soft ideals" the codimension conjecture always holds:

\begin{corollary}\label{C: soft-edged or soft-complemented}
If $I$ is a soft-edged or soft-complemented ideal, then
\[
\di \frac{I}{[I,B(H)]} \quad \text{is} \quad
\begin{cases}
1 \quad \text{or uncountable} \quad    &\text{if $\omega \notin \Sigma(I)$} \\
0 \quad \text{or uncountable} \quad    &\text{if $\omega \in \Sigma(I)$.}
\end{cases}
\]
\end{corollary}

\noindent In particular, $\di \frac{(\omega)}{[(\omega),B(H)]}$ is uncountable since $(\omega)$ is soft-complemented as all principal ideals are, but is not am-stable since  $\omega\not\asymp\omega_a$ is not regular. \\

Further consequences and extensions are:

\begin{corollary}\label{C:am/am-infty uncountable dimension}
Let $I \ne \{0\}$ be an ideal.
Then $\frac{I}{[I,B(H)]}$ has uncountable dimension if any of the following conditions hold.
\item[(i)] $I \not\supset st^a(\mathscr L_1)$ and
$I \subset J$ but $I\not\subset st_{a_\infty}(J)$ for some soft-complemented ideal $J$.
\item[(ii)] $I\not\subset st_{a_\infty}(\mathscr L_1)$ and
$J \subset I$ but $st^a(J)\not\subset I$ for some soft-edged ideal $J$.
\end{corollary}

One case of interest is when $I$ or $J$ are the trace class $\mathscr L_1$, which is both soft-edged and soft-complemented.

\begin{corollary}\label{C:dimensions of 3 quotients} Let $I$ be a nonzero ideal. Then
\item[(i)]
\[
\di \frac{I+[\mathscr L_1,B(H)]}{[\mathscr L_1,B(H)]} \quad \text{is} \quad
\begin{cases}
1     &\text{if $I \subset ~_{a_\infty}(\mathscr L_1)$} \\
\text{uncountable}    &\text{if $I \not\subset ~_{a_\infty}(\mathscr L_1)$}
\end{cases}
\]
\item[(ii)]
\[
\di \frac{\mathscr L_1+[I,B(H)]}{[I,B(H)]} \quad \text{is} \quad
\begin{cases}
0     &\text{if $\omega \in \Sigma(I)$} \\
1     &\text{if $\omega \in \Sigma(\sc I) \setminus \Sigma(I)$} \\
\text{uncountable}    &\text{if $\omega \notin \Sigma(\sc I)$}
\end{cases}
\]
\item[(iii)] If $\omega \notin \Sigma(I)$ then
\[
\di \frac{I}{I \cap \mathscr L_1+[I,B(H)]} \quad \text{is} \quad
\begin{cases}
0     &\text{if $I \subset \mathscr L_1$} \\
\text{uncountable} \quad    &\text{if $I \not\subset \mathscr L_1$}
\end{cases}.
\]
In particular, if $\mathscr L_1 \not\subset I$, then there are uncountably many linearly independent extensions of $Tr$ from $\mathscr L_1$ to $I$.
\end{corollary}

\begin{remark}\label{R: Dixmier proved}
Dixmier proved in \cite{jD66} that the am-closure of a principal ideal, $(\eta)^- = ~_a(\eta_a)$,
for which $(\eta) \subset \se (\eta)_a$ (i.e., $\eta = o(\eta_a)$, which is equivalent to $\frac{(\eta_a)_{2n}}{(\eta_a)_n} \rightarrow \frac{1}{2}$) supports a positive singular trace.
In \cite[Section 5.27 Remark 1]{DFWW} it was noted that Dixmier's construction can be used to show that 
$\di \frac{(\eta)^-}{\cl [(\eta)^-,B(H)]} = \infty$
where $\cl$ denotes the closure in the principal ideal norm.
Corollary \ref{C: soft-edged or soft-complemented} shows that $\di \frac{(\eta)^-}{[(\eta)^-,B(H)]}$ is uncountable
follows from the weaker hypothesis $\eta_a \ne O(\eta)$,
i.e., $\eta$ is not regular.
\end{remark}

Another case where we can prove that $[I,B(H)]$ has infinite codimension (we do not know if uncountable) is the following. 
Its proof is unrelated to the proof  of Theorem \ref{T: BASICS: seJ notin F+[I,B(H)] implies uncountable dimension-(J+[I,B(H)])/(F+[I,B(H)])} but is more combinatorial in nature.

\bC{C:otherinfcodim}\cite[Corollary 7.9]{vKgW04-Traces}
Let $I$ be an ideal for which there exists a nonsummable sequence $\xi \in \Sigma(I)$ and
a monotone sequence of indices $\{p_k\}$ so that for every $\chi \in \Sigma(I)$ there is an associated
$m \in \mathbb N$ for which $(\frac{\chi}{\xi_a})_{mp_k} \rightarrow 0$. Then $\di \frac{I}{[I,B(H)]} = \infty$.
\eC

This result and the unusual property for the ideal $I$ originates from

\begin{theorem}\label{T: false converse for T: se or sc inf codim implies}\cite[Theorem 7.8]{vKgW04-Traces}
For every am-stable principal ideal $J \ne \{0\}$ there is an ideal $I$ with 
$\se J \subset I \subset J$ for which $\di \frac{I}{[I,B(H)]} = \infty$, yet $\se I$ and hence $\sc I$ are am-stable.
\end{theorem}

\section{\leftline{\bf B(H) lattices }}\label{S: Density}

\subsection{\leftline{Lattice density}}\label{S:1} \quad
Many aspects of operator theory implicitly (and sometimes explicitly) involve the structure of the lattice $\mathscr L$ of all operator ideals and of its distinguished sublattices, like the lattice $\mathscr PL$ of all principal ideals. For instance, one basic question due to Brown, Pearcy and Salinas \cite{BPS71} was whether or not the ideal of compact operators $K(H)$ was the sum of two proper ideals. This was settled affirmatively in \cite{aBgW78} for any proper ideal strictly larger than $F$ (the ideal of finite rank operators) using the continuum hypothesis and the techniques employed led to the set-theoretic concept of groupwise density that has proved useful in point-set topology and abelian group theory \cite{aB89}, \cite{aB90}, \cite{hM01}. Other work by Salinas (e.g., \cite{nS74}) studied increasing or decreasing nests of special classes of symmetric norm ideals (Banach ideals). Furthermore, many of the questions on the structure of arithmetic mean ideals depend on properties of operator ideal lattices, particularly various lattices of principal ideals.

Operator ideals form a lattice $\mathscr{L}$ with inclusion as partial order and intersection as meet and sum as join.
Since, for any $\xi,\eta \in \co*$, $(\xi) \cap (\eta) = (\min(\xi,\,\eta) )$ and $(\xi) + (\eta) = (\xi+\eta )$,
the collection of all principal ideals forms a sublattice which we denote by $\mathscr{PL}$.
Since the intersections and sums of two am-stable ideals are easily seen also to be am-stable,
the collection of all am-stable principal nonzero ideals forms a lattice that we denote by $\mathscr{SPL}$.
And likewise for the am-$\infty$ case which lattice we denote by $\mathscr{S_{\infty}PL}$. 
Similarly, the collection of all nonzero principal ideals having a generator sequence that satisfies the \D* forms a lattice denoted by $\Delta_{1/2}\mathscr{PL}$.

Lattices may have gaps: a nested pair of ideals in the lattice, $I \subsetneq J$, between which there is no ideal in the lattice. 
In fact, given any principal ideal $(\xi)$, an easy maximality argument provides us with an ideal $I \subsetneq  (\xi)$ that forms a gap under $(\xi)$. 
On the other hand, there cannot be any gap above $(\xi)$, i.e., $(\xi) \subsetneq J$. 
We describe this fact by saying that  $\mathscr{PL}$ is \textit{upper dense} but not \textit{lower dense} in $\mathscr{L}$ according to the following definition. 

\begin{definition}\label{L: density, strong density and gaps}
 For two nested lattices $\mathscr L' \subset \mathscr
L''$, $\mathscr L'$ is upper dense (resp.,
lower dense) in $\mathscr L''$ provided that
between every pair of ideals $I \in \mathscr L'$ and  $J\in \mathscr
L''$ with $I\subsetneq J$ (resp., $J\subsetneq I$) lies another ideal in the smaller lattice $\mathscr L'$.
\end{definition}

Many questions on operator ideals can be reduced to questions on lattice density. Together with the upper density of $\mathscr {PL}$ in $\mathscr {L}$, the following three theorems that are the technical core of \cite{vKgW04-Density}, 
provided the tools we need for applications to arithmetic mean ideals. 

\begin{theorem}\label{T: Delta1/2-PL strong density and strong gaps}\cite[Theorem 3.9]{vKgW04-Density}
$\Delta_{1/2} \mathscr {PL}$ is upper and lower dense in $\mathscr {PL}$.
\end{theorem}

\begin{theorem}\label{T: SPL density in PL}\cite[Theorem 3.10]{vKgW04-Density}
$\mathscr {SPL}$ is upper and lower dense in $\mathscr {PL}$.
\end{theorem}

\begin{theorem} \label{T: SinftyPL density}\cite[Theorem 3.13]{vKgW04-Density}
$\mathscr {S_{\infty}PL}$ is upper and lower dense in $\mathscr {PL}$.
\end{theorem}

A key property used in the proofs is the following Potter-type condition characterizing a regular sequence $\xi $: $\xi_n\ge C(\frac{m}{n})^{p_o}\xi_m$ for some $0< C \le 1, 0<p_o<1$ and all $n\ge m$ 
\cite[Theorem 3.10 and Remark 3.11]{DFWW}. 
A similar key role is taken by the characterization of an $\infty$-regular sequence $\xi $: $\xi_n\le C(\frac{m}{n})^{p_o}\xi_m$ for some $C \ge 1, p_o>1$ and all $n\ge m$ given \linebreak
\cite[Theorem 4.12]{vKgW04-Traces}.

Since $\mathscr {PL}$ is upper but not lower dense in $\mathscr {L}$, immediate consequences of Theorems \ref{T: Delta1/2-PL strong density and strong gaps}, \ref{T: SPL density in PL}, and \ref{T: SinftyPL density} are:

\begin{corollary} \label{C:density in L}
\item[(i)] $\mathscr {PL}$, $\Delta_{1/2} \mathscr {PL}$, $\mathscr{SPL}$, and $\mathscr {S_{\infty}PL}$
are upper dense in $\mathscr {L}$ but they all have lower gaps in
$\mathscr {L}$.
\item[(ii)] $\mathscr {PL}$, $\Delta_{1/2} \mathscr {PL}$, $\mathscr
{SPL}$, and $\mathscr {S_{\infty}PL}$
   have no gaps.
\end{corollary}

Other consequences of the same theorems and the Potter-type
conditions already mentioned are:

\begin{corollary}\label{C:inclusions}
 \item[(i)] Every principal ideal is contained in a principal ideal
from $\mathscr {SPL}$,
and if it is strictly larger than $F$ it contains a principal ideal
in $\mathscr {S_{\infty}PL}$ that is strictly
  larger than $F$.
\item[(ii)] A principal ideal $(\eta)$ contains a principal ideal in
$\Delta_{1/2} \mathscr {PL}$ if and only
  if $\omega^p = O(\eta )$ for some $p>0$; it contains a principal
ideal in $\mathscr {SPL}$  if and only if
$\omega^p = O(\eta )$ for some $0 < p < 1$; and it is contained in  a
principal ideal in $\mathscr
{S_{\infty}PL}$ if and only if $\eta = O(\omega^p)$ for some $ p > 1$.
\end{corollary}

We can reformulate some of these results  in terms of traces. 
By Corollary \ref {C:inclusions}(i), every principal ideal $I$ is contained in a principal ideal $J$ with no nonzero traces; 
thus no nonzero trace on $I$ can extend to a trace on $J$.
Alternatively, by Corollary \ref{C:inclusions}(ii), 
if a principal ideal $I$ is sufficiently large to contain $(\omega^p)$ for some $0<p<1$, then $(\omega^p)= [(\omega^p), B(H)] \subset [I, B(H)]$ and hence every trace on $I$ must necessarily vanish on $(\omega^p)$.
More generally, by employing repeatedly Theorem \ref{T: SPL density in PL} and \cite[Proposition 3.17]{vKgW04-Density}, every am-stable principal ideal is the first ideal of an increasing (or decreasing) countable chain of principal ideals where every odd numbered ideal is am-stable and hence has no nonzero trace, and every even numbered ideal supports infinitely many nonzero traces, 
that is, its space of traces is infinite-dimensional
(\cite[Theorem 7.5]{vKgW04-Traces} and \cite[Proposition 4.6(i)]{vKgW04-Soft}).

The am-$\infty$ case is similar.
\subsection{\leftline {Unions of principal ideals} }\label{S:4} \quad
Basic questions on intersections and unions of ideals from a certain lattice play a natural role in the subject.
Salinas and others investigated intersections and unions of ideals related to various classes of mainly Banach ideals (e.g., see \cite {BPS71}, \cite{nS74}).
We used Corollary \ref{C:inclusions} and density 
Theorems \ref{T: Delta1/2-PL strong density and strong gaps}-\ref{T: SinftyPL density} to determine intersections and unions of various classes of principal ideals.
Then we investigated questions on representing ideals as unions of chains of principal or countably generated ideals.\\

An immediate consequence of Corollary \ref {C:inclusions} is:

\begin{corollary}\label{L:intersections} \quad\\
(i) $\underset{p >0} \bigcap\,(\omega^p) = \underset{L \in \Delta_{1/2} \mathscr {PL}}\bigcap L$ 
\qquad \qquad (ii) $\underset{0 < p < 1} \bigcap\,(\omega^p) = \underset{L \in \mathscr {SPL}}\bigcap L$ 
\qquad \qquad (iii) $\underset{p > 1} \bigcup\,(\omega^p) =  \underset{L \in \mathscr {S_{\infty}PL}}\bigcup L$
\end{corollary}

In the case of Banach ideals or their powers we can obtain more:
\begin{proposition}\label{P:complete}\cite[Proposition 4.3]{vKgW04-Density}
\item(i) Let $I=J^p$ for some $p>0$ and some Banach ideal $J$.
Then $I = \bigcup \{L \in \Delta_{1/2}\mathscr {PL} \mid L \subset I\}$.

\item(ii) For a Banach ideal $I$, $I$ is am-stable if and only if $I = \bigcup \{L \in \mathscr {SPL} \mid L \subset I\}$.

\item(iii) For a Banach ideal $I$, $I$ is am-$\infty$ stable if and only if $I = \bigcup \{L \in \mathscr {S_{\infty}PL} \mid L \subset I\}$.
\end{proposition}

While every ideal is the union of principal ideals,
a different and natural question is: which ideals are the union of an increasing chain of principal ideals?
A partial answer is given by the following proposition.

\begin{proposition}\label{P: nested unions and the continuum hypothesis} \cite[Proposition 4.4]{vKgW04-Density} \\
Assuming the continuum hypothesis,
the following hold.
\item[(i)] All ideals are unions of increasing chains of countably
generated ideals.
\item[(ii)] An ideal $I$ is the union of an increasing chain of
principal ideals if for every countably generated ideal
  $J\subset I$ there is a principal ideal $L$ with $J\subset L\subset I$.
\end{proposition}

\cite[Corollary 4.5]{vKgW04-Density} shows that Banach ideals and more generally $e$-complete ideals,
i.e., powers of Banach ideals (see \cite[Section 4.6]{DFWW}), satisfy the condition in part (ii) of the above proposition.

Neither the $e$-completeness condition in the corollary nor the
condition in part (ii) of the proposition are necessary and we do not know whether or not the continuum hypothesis or
something stronger than the usual axioms of set theory
is required for Proposition \ref {P: nested unions and the continuum hypothesis}.

Clearly, an ideal is the union of a \textit{countable} increasing chain of principal ideals if and only if it is countably generated. But which ideals cannot be  the union $I=\cup_{k=1}^\infty J_k$ of a countable increasing chain of ideals  $J_k \subsetneq I $? Salinas showed in  \cite[Theorem 7.4] {nS74} that every minimal Banach ideal $I$, i.e., an ideal  $I=\mathfrak{S}_\phi^{(o)}$ given by a symmetric norming function $\phi$, has this property. Maximal Banach ideals may fail this property: e.g., consider any principal am-stable ideal.  Conversely,   maximal Banach ideals, i.e., ideals $\mathfrak{S}_\phi$ given by a symmetric norming function $\phi$, cannot be the intersection of a countable decreasing chain  of \textit{norm} ideals that properly contain them \cite[Corollary  7.8] {nS74}.  To the best of our knowledge not much more progress has been achieved since then on this line of questions.

\subsection{\leftline {Applications to first order arithmetic mean ideals} }\label{S:5} \quad
While in general the first order arithmetic mean ideals can be am-stable without the ideal $I$ being am-stable, e.g., the examples mentioned at the end of Section \ref {S:5} (and likewise for am-$\infty$ ideals), principal ideals are more ``rigid."

First order arithmetic mean ideals are investigated in \cite{DFWW} and \cite {vKgW02}-\cite{vKgW04-Soft} and the density theorems of Section 3 provide further information for the principal ideal cases.

\begin{theorem}\label{T: I am-stable iff}\cite[Theorem 5.1]{vKgW04-Density}\\
Let $I$ be a nonzero principal ideal. 
Then am-stability of the following ideals are equivalent: $I,I^{oo},I_a$ and $I^-$, and also $_aI,I^o$ and $I_-$ provided the latter three are nonzero.
\end{theorem}

\noindent Notice that for any ideal $I$ that does not contain
$(\omega)$ (resp., $\mathscr L_1$), $_aI=I^o = \{0\}$
(resp., $I_-= \{0\}$). Since $\{0\}$ is am-stable while in both cases
$I$ is not,
  the conditions in the last corollary that $I_a$, $I^o$ and $I_-$
are nonzero is essential. Some of the implications in Theorem \ref {T: I am-stable iff}\ depend on Proposition \ref {P:open'} below. \\

There are two main differences in the am-$\infty$ case.
The first is that if $I \neq \{0\}$, then unlike the am-case, $_{a_{\infty}}I$ is never zero, and in lieu of the two chains of inclusions, we have
$_{a_{\infty}}I \subset I^{o\infty} \subset I$, $_{a_{\infty}}I \subset I_{-\infty} \subset I$, $I\cap \mathscr L_1 \subset I^{-\infty} \subset I_{a_{\infty}}$, and $I \cap \se (\omega) \subset I^{oo\infty} \subset I_{a_{\infty}}$. \linebreak
(See \cite[Proposition 4.8(i)-(i$'$)]{vKgW04-Traces} and remarks following both \cite[Corollary 3.8 and Proposition 3.14]{vKgW04-Soft}.)
The second is that if $I$ is principal,
then $I^{o\infty}$, $I^{oo\infty}$ and $I_{a_{\infty}}$ are either principal or they coincide with $\se (\omega)$
\cite[Lemma 4.7]{vKgW04-Traces}, \cite[Lemma 3.9, Lemma 3.16]{vKgW04-Soft}.
Consequently, as $\se (\omega)$ is not am-stable, when $I$ is principal, if any of the ideals $I^{o\infty}$, $I^{oo\infty}$ or $I_{a_{\infty}}$ is am-stable, then it is principal as well.

\begin{theorem}\label{T: I am-inf stable iff}\cite[Theorem 5.2]{vKgW04-Density}
Let $I \neq \{0\}$ be a principal ideal. Then am-$\infty$ stability
of the following ideals is equivalent:
$I$, $_{a_{\infty}}I$, $I^{o{\infty}}$, $I_{-{\infty}}$,
$I^{oo\infty}$, $I_{a_{\infty}}$, and $I^{-\infty}$.
\end{theorem}

The fact that a sequence $\xi \in \co* $ is regular if and only if $\xi_a$ is regular has been proved in a number of different ways 
(see \cite[Remark 3.11]{DFWW} for a discussion).
The analogous result for $\infty$-regularity \cite[Theorem 4.12]{vKgW04-Traces} is quite hard and was one of the main results in that paper.
Both facts are also immediate consequences of Theorems \ref{T: I am-stable iff} and \ref{T: I am-inf stable iff}.\\

It is useful in the study of first order cancellations to note that Theorem \ref{T: I am-stable iff}
does not extend beyond principal ideals. In fact, in \cite[Examples 5.4-5.5] {vKgW04-Density}  we find a countably generated ideal $N$ which is not am-stable but for which
$_aN =N^o =N_-= st^a((\omega))$ is am-stable and a countably generated ideal $L$ which is not am-stable but for which $L_a =L^- =L^{oo} = st^a((\omega))$ is am-stable.\\

\noindent \textbf{Question 6.}
Find am-$\infty$ analogs for $N$ and $L$, or prove otherwise.\\

\subsection{\leftline {First order cancellation properties}} \label{S:6} \quad
Which ideals $I$ have the  \textit{cancellation  property} that $J_a\subset I_a \Rightarrow J\subset I$? It is immediate to verify that the answer is: all am-closed ideals. Similarly, the first order cancellations $_aI\subset \,_aJ\Rightarrow I\subset J$ characterize the am-open ideals, and analogous results hold for the am-$\infty$ case \cite[Lemma 5.13]{vKgW04-Density}.
 
 As a consequence, if $I$  is a countably generated ideal that is not am-stable, it is not am-closed by Theorem  \ref  {T:coun gen closed} and hence there is some ideal $J\not\subset I$ for which $J_a\subset I_a$. But then $L:=I+J\supsetneq I$ with $L_a=I_a+J_a=I_a$. In other words, equality cancellation does not hold when $I$ is not am-stable.
 
 At an early stage of this project, Ken Davidson and the second named author found a direct constructive proof of this fact for the case of principal ideals or, equivalently, that there is a principal ideal $J\supsetneq I$ with $J^-=I^-$. 
The same result was obtained earlier by Allen and Shen \cite{AS78} by different methods. 
 
Notice that if $I$  is principal and not am-stable, then it may be impossible to find an ideal $J \subsetneq I$ such that $J_a=I_a$. 
A counterexample is constructed in \cite[Example 6.3]{vKgW04-Density}. The same positive and negative results hold in  the am-$\infty$ case.
 
The question of whether for a countably generated ideal $I$ there exists an ideal $J\supsetneq I$ with 
with $_aJ=\,_aI$ is considerably harder even when we assume that $I$ is principal and it requires a direct constructive proof. We first notice that  the condition $_aJ=\,_aI$ is equivalent to $J^o=I^o$. The answer for principal ideals is given by;
\begin{proposition}\label{P:open'}\cite [Proposition 6.5]{vKgW04-Density}
\item[(i)]
Let  $(\rho)\ne\{0\}$ be a principal ideal that is not am-stable.
Then there are two principal ideals $(\eta^{(1)})$ and $(\eta^{(2)})$
with $(\eta^{(1)})$ possibly zero, such that
\[
(\eta^{(1)}) \subsetneq (\rho) \subsetneq (\eta^{(2)}),\quad
(\eta^{(1)})^{oo} = (\rho)^o, \quad \text{and} \quad
(\rho)^{oo}=(\eta^{(2)})^o.
\]
\item[(ii)] Let  $(\rho)\ne\{0\}$ be a principal ideal that is not am-$\infty$ stable.
Then there are two principal ideals $(\eta^{(1)})$ and $(\eta^{(2)})$ such that
\[
(\eta^{(1)}) \subsetneq (\rho) \subsetneq (\eta^{(2)}), \quad
(\eta^{(1)})^{oo\infty} = (\rho)^{o\infty}, \quad \text{and} \quad
(\rho)^{oo\infty}=(\eta^{(2)})^{o\infty}.
\]
\end{proposition}
This proposition was the key for proving some of the implications in Theorems \ref {T: I am-stable iff} and \ref {T: I am-inf stable iff} presented above and also the following theorem.

\begin{theorem}\label{T: cancellation}\cite[Theorem 6.7]{vKgW04-Density}
\item[(A)] Let $I\ne \{0\}$ be a principal ideal and let $J$ be an
arbitrary ideal. Then the following are equivalent.
\item[(i)] $I$ is am-stable
\item[(ii)] $_aJ = \,_aI$  (or, equivalently, $ J^o = I^o$) implies $J = I$
\item[(iii)] $ J_a = I_a$  (or, equivalently, $ J^- = I^-$) implies $J = I$

\item[(B)] Let $I\ne \{0\}$ be a principal ideal and let $J$ be an
arbitrary ideal. Then the following are equivalent.
\item[(i)] $I$ is am-$\infty$ stable
\item[(ii)] $_{a_\infty}J =\, _{a_\infty}I$ (or, equivalently, $
J^{o\infty} = I^{o\infty}$) implies $J = I$.
\item[(iii)] $J_{a_\infty} =I_{a_\infty}$  (or, equivalently, $
J^{-\infty} = I^{-\infty}$) implies $J = I$
\end{theorem}
\noindent  Similar cancellation properties hold for the identities  $J^{oo}=I^{oo}$,  $J_-=I_-$ and for the am-$\infty$ case.\\

The examples mentioned at the end of Section \ref {S:5} show that the first order equality cancellations of Theorem \ref{T: cancellation} can fail for countably generated am-stable ideals.
Indeed, $_aN=N_-$ is am-stable yet 
\[
_aN=\,_{a^2}N \not \Rightarrow N=\,_aN \text{ and } N_-=\,(N_-)_- \not \Rightarrow N=\,N_-. 
\]
Similarly, $L_a = L^{oo}$ is am-stable yet $L_a=L_{a^2}$ does not imply $L=L_a$ and $ L^{oo}=( L^{oo}) ^{oo}$ does not imply $L=L^{oo}$.\\

First order inclusion  cancellations occur naturally in the study of operator ideals. We have seen that two of the inclusion cancellations characterize am-open and am-closed ideals (and similarly for the am-$\infty$ case.  Another  cancellation question, ``for which ideals $I $ does $J_a\supset I_a$ imply $J \supset I$?"  has an interesting origin and an answer that we also found interesting.

Its origin is the still open question from \cite[Section 7]{DFWW}:
If the class $[J,B(H)]_1$ of single commutators of operators in $J$ with operators in $B(H)$,
contains a finite rank operator with nonzero trace, must $J\supset (\sqrt \omega)$?

It was shown in \cite[Theorem 7.3]{DFWW} only that then  $J^- \supset (\sqrt \omega)= (\sqrt \omega)_a$, so it was natural to ask whether this result was enough, namely, whether  $(\sqrt \omega)$ satisfied this cancellation property.  To obtain the answer (negative), we define for any ideal $I$ the   ``greatest lower bound," $\widehat I$ for this cancellation, namely:

\begin{definition}\label{D:Gg} 
$\widehat{I}:=\bigcap\{J \mid J_a \supset I_a\}$
\end{definition}
\noindent We can explicitly identify $\widehat{I}$ when $I$ is principal and show that it is itself principal.
\begin{definition}\label{D: 1}
For $\xi \in \text{c}_\text{o}^* \setminus \ell^1$,
let $\nu(\xi)_{n} := \min \{k \in \mathbb N \mid
\sum_{i=1}^{k}\xi_{i} \geq n\xi_1 \}$ and define
$\widehat{\xi} :=\, <(\xi_{a})_{\nu(\xi)_{n}} >$.
\end{definition}

\noindent Then $\widehat{\xi}$ is asymptotic to $<\frac{n}{\nu(\xi)_n}\xi_1>$ and moreover

\begin{theorem}\label{T: 4}\cite [Theorem 6.14]{vKgW04-Density}
Let $\xi \in \co* \setminus \ell^1$.
\item[(i)] If $\eta_a \geq \xi_a$ and $\eta \in \co*$, then for every
$\varepsilon > 0$,
$\eta_{n} \geq (1 - \varepsilon) \,\widehat{\xi}_{n}$
for $n$ sufficiently large.
\item[(ii)]
For every $\rho \in \co*\setminus \Sigma((\widehat \xi))$,
there is an $\eta \in \co*$ for which $\eta_a \geq \xi_a$ and
$\rho \notin \Sigma((\eta))$.
\end{theorem}

As a consequence we identify $\widehat I$ for all  principal ideal $I$ and in particular, for all ideals $(\omega^p)$.

\begin{corollary}\label{C: 5} \cite[Corollary 6.15]{vKgW04-Density}
Let $0 \ne \xi \in \co*$, then
$\widehat{(\xi)} =
\begin{cases}
F &\text{\quad for~$\xi \in \ell^1$}\\
(\widehat\xi) &\text{\quad for~$\xi\notin \ell^1$}
\end{cases}$
\end{corollary}

\begin{corollary}\label{C: 6} \cite[Corollary 6.16]{vKgW04-Density}
\item[(i)]
\[ \widehat{(\omega)^p} =
\begin{cases}
F &\text{for $p>1$}\\
(<\frac{n}{e^n}>) &\text{for $p=1$}\\
(\omega)^{p'} &\text{for $0<p<1$, where $\frac{1}{p} -\frac{1}{p'}=1$}
\end{cases}
\]
\item[(ii)] $\widehat{(\omega)^p}$ is am-stable if and only if
$0<p<\frac{1}{2}$.
\end{corollary}

So, in particular $\widehat{(\omega)^{1/2}} = (\omega)$: in other words, knowing that $J^-\supset  (\omega)^{1/2}$ only guarantees that $J\supset (\omega)$.

Notice that by definition, $\widehat{(\xi)}\subset (\xi)$, but as the
case of $(\omega)^{1/2}$ illustrates, the inclusion can be proper
even when $\xi$ is regular, i.e., when $(\xi)$ is
am-stable. As the following proposition shows, the inclusion is
certainly proper when $(\xi) \neq F$ is not am-stable.

\begin{proposition}\label{P:Gg is stable}\cite [Proposition 6.17] {vKgW04-Density}

 If $F\neq (\xi) = \widehat{(\xi)}$ with $\{0\}\neq\xi \in \co* $, then $(\xi) =(\xi)_a$.
\end{proposition}
 
 There are principal ideals $ I=\widehat{I}$ with $I\neq F$.

\begin{example}\label{E:7} \cite[Example 6.18]{vKgW04-Density}
 Let $\xi_j = \frac{1}{k}$ for $j\in
((k-1)!, k!]$. Then $(\xi) = \widehat{(\xi)}$.
\end{example}
 
Thus the condition $\xi\asymp \widehat{\xi}$ is strictly stronger than regularity. \\

While the  definition of $\widehat{I}$ can be extended naturally to the am-$\infty$ case by setting 
$\widehat{I}^\infty:=\,\bigcap\{J \mid J_{a_\infty} \supset I_{a_\infty}\}$, we have not analyzed its properties.\\

\noindent \textbf{Question 7.}
If $I=(\xi)$, is $\widehat{I}^\infty$ principal? If yes, is there an expression for its generator similar to the one in Definition \ref {D: 1} and with similar properties?\\

\section{\leftline{\bf Second order arithmetic means in operator ideals}}\label{S: 2nd Order}
Second order arithmetic mean cancellations are considerably more complex even for principal ideals and are the focus of \cite{vKgW04-2nd Order}.
The questions we addressed there are: which conditions on an ideal $I$ guarantee that the following second order arithmetic mean inclusion cancellations and equality cancellation hold for arbitrary $J$?
\begin{itemize}
\item[(i)] $J_{a^2} \subset I_{a^2} \Rightarrow J_a \subset I_a$
\item[(ii)] $J_{a^2} \supset I_{a^2} \Rightarrow J_a \supset I_a$
\item[(iii)] $J_{a^2} = I_{a^2} \Rightarrow J_a = I_a$ \\
\end{itemize}

The first natural ``test" question, posed by  M. Wodzicki arising from work in \cite{DFWW},
is whether equality cancellation (iii) holds automatically for all pairs of principal ideals.
Reformulated in terms of the s-number sequences $\xi$ and $\eta$ of the generators of the two principal ideals, the question asks whether
the equivalence $\xi_{a^2}\asymp \eta_{a^2}$ of the sequences of the second order arithmetic means always implies the equivalence $\xi_a \asymp \eta_a$ of the first order means.

The answer to this question is negative and is one of the main results of \cite{vKgW04-2nd Order}.

\begin{example}\label{E: Example 6}\cite[Example 4.3]{vKgW04-2nd Order}
There exists principal ideals $J \subset I $ for which $J_{a^2}= I_{a^2}$ but $J_a\ne I_a$.

\end{example}
\begin{proof}[Construction] Set $\delta_k=e^{-\sum_{p=1}^{k-1}p^2}$ for $k>1$, $\xi_1=\eta_1 = 1$,  and define
 \[
\xi_j:= e^{-k^2}\delta_k ~\text{ for }~ m_k < j \le m_{k+1}
\]
and
\[
\eta_j :=
\begin{cases}
ke^{-k^2}\delta_k     &\text {for  $m_k < j \le [e^{k^2}m_k]$} \\
e^{-k^2}\delta_k    &\text{for $[e^{k^2}m_k] < j \le  m_{k+1}$.}
\end{cases}
\]
where  $m_{k+1}> [e^{k^2}m_k]$ is chosen sufficiently large so to insure that $(\eta_{a^2})_{m_{k+1}} \le
(1+ \frac{1}{k})e^{-k^2}\delta_k$.  Then
\[
\frac{(\xi_{a^2})_{n_k}}{\delta_k} \sim \frac{(\eta_{a^2})_{n_k}}{\delta_k}\sim k^2e^{-k^2}
\]
while
\[
\frac{(\xi_a)_{n_k}}{\delta_k}  \sim 2e^{-k^2}\quad \text {and}\quad \frac{(\eta_a)_{n_k}}{\delta_k}  \sim
k e^{-k^2}.
\]
Let $J:=(\xi)$ and $I:=(\eta)$, then $J\subset I$, $J_{a^2}= I_{a^2}$ but $J_a\neq I_a$.
\end{proof}

\noindent The intuition behind the construction of this example led to the notion
of the \textit{ratio of regularity} sequence
$r(\xi): = \frac{\xi_a}{\xi}$ for a nonincreasing sequence $\xi \in \text{c}_{\text{o}}$ (see \cite[Section 3]{vKgW04-2nd Order}) and then, indirectly, to the other results in \cite{vKgW04-2nd Order}.

Notice that in general the ratio of regularity $\frac{\xi_a}{\xi}$ has downward variations bounded by the inequality $r_{n+1} \geq  \frac{n}{n+1}r_n+ \frac{1}{n+1} $, but it can vary abruptly upwards because the ratio $\frac{r(\xi)_{n+1}}{r(\xi)_{n}} \ge \frac{n}{n+1}\frac{\xi_n}{\xi_{n+1}}$
can be arbitrarily large. This is not the case, however, if $\xi$ itself is an average. Indeed
$ r(\eta_a)_{n+1}  < (1+\frac{1}{n})r(\eta_a)_n$ for all $n$ and all $\eta$. Moreover,
\[
1\le r(\eta_a)_n \le \log n\quad \text{ for all $n >1$ and all }\eta.
\]
The two ``extremal cases" of  the ratio of regularity of an average are thus when it is bounded (i.e.,
$ r(\eta_a)\asymp 1$) and when it is ``maximal," (i.e.  $r(\eta_a) \asymp log$). The first case occurs precisely when $\eta_a$ is regular (if and only if $\eta$ is regular). The second case is equivalent  to what we call the exponential $\Delta_2$-condition
$\underset{m}{\sup}\,\frac{m^2(\eta_a)_{m^2}}{m(\eta_a)_m} < \infty$ \cite[Proposition 3.11]{vKgW04-2nd Order}.

Somewhat to our surprise, we found that either extremal cases for $\eta$ are sufficient, in the principal ideal case, for the equality cancellation in (iii) to hold:

\bT {T: Theorem 4.8}\cite [Theorem 4.8]{vKgW04-2nd Order}
\item[(i)]  Let $I = (\eta)$ be principal ideal, let $r(\eta_a)\asymp
\log$  or $\eta \asymp \eta_a$,  and let $J$ be an arbitrary ideal.
Then  $J_{a^2} = I_{a^2} \Rightarrow J_a = I_a$.
\item[(ii)] Let $I$ be an ideal such that for every $\xi \in \Sigma(I)$ there is some $\eta \in \Sigma(I)$, $\xi\le \eta$, for which $r(\eta_a) \asymp \log$ or $\eta \asymp \widehat \eta$  and let $J$ be an arbitrary ideal. Then  $J_{a^2} = I_{a^2} \Rightarrow J_a = I_a$.
\eT
As mentioned at the end of Section \ref {S:6} the condition $\eta \asymp \widehat \eta$ is strictly stronger than  $\eta \asymp \eta_a$. We do not know if the weaker condition $\eta \asymp  \eta_a$ might also suffice in part (ii) of Theorem \ref{T: Theorem 4.8}. Some partial negative evidence comes from the fact that for general ideals (even countably generated ideals), arithmetic mean stability is not sufficient for the first order equality cancellation to hold, and the counterexample is the ideal $L$ presented after Theorem \ref {T: cancellation}. However, $L$ itself is not the arithmetic mean of an ideal, i.e., it is not am-open, so that it does not provide a counterexample to second order equality cancellation.\\

\noindent \textbf{Question 8.}
\begin{enumerate}
 \item[(i)] Would the condition  $\eta \asymp  \eta_a$ be sufficient in part (ii) of Theorem \ref{T: Theorem 4.8}?
 \item[(ii)] What is a necessary condition for the principal ideal case?
 \end{enumerate}

While we do not know if these two conditions, $\eta$ regular or $r(\eta_a)\asymp \log$, are necessary for the equality cancellation in (iii) to hold, we know that they are too strong for the inclusion cancellation in (i). Indeed a weaker  sufficient condition for (i) is that $r(\eta_a)$ is equivalent to a monotone sequence.

\bT{T: Theorem 4.5}\cite[Theorem 4.5] {vKgW04-2nd Order}
\item[(i)] Let $I = (\eta)$ be principal ideal, let $r(\eta_a)$ be equivalent to a monotone sequence and let $J$ be an arbitrary ideal.
Then $J_{a^2} \subset I_{a^2} \Rightarrow J_a \subset I_a$.
\item[(ii)]  Let $I$ be an ideal such that for every $\xi \in \Sigma(I)$ there is some $\eta \in \Sigma(I)$, $\xi\le \eta$,  for which $r(\eta_a)$ is equivalent to a monotone  sequence and let $J$ be an arbitrary ideal. Then  $J_{a^2} \subset I_{a^2} \Rightarrow J_a \subset I_a$.
\eT

On the other hand, the two conditions, $\eta$ regular or $r(\eta_a)\asymp \log$, are too weak for the inclusion cancellation in (ii).
Indeed, the principal ideal generated by the regular sequence $\omega^{1/2}$ does not satisfy the inclusion cancellation in (ii).

\begin{example}\label{e: omega1/2}  \cite[Example 4.7]{vKgW04-2nd Order}
There exists a principal ideal $(\xi)$ for which $(\omega^{1/2})_{a^2}\subset (\xi)_{a^2}$ but for which $(\omega^{1/2}) _{a}\not\subset (\xi)_{a}$.
\end{example}

A sufficient condition for the inclusion cancellation in (ii) is:

\bT{T: Theorem 4.6}\cite[Theorem 4.6]{vKgW04-2nd Order}
\item[(i)] Let $I = (\eta)$ be a principal ideal, let $r(\eta_a)\asymp \log$ or $\eta \asymp \widehat \eta$, and let $J$ be an arbitrary ideal.
Then $J_{a^2} \supset I_{a^2} \Rightarrow J_a \supset I_a$.
\item[(ii)] Let $I$ be an ideal such that for every $\xi \in \Sigma(I)$ there is some $\eta \in \Sigma(I)$, $\xi\le \eta$, for which $r(\eta_a) \asymp \log$ or $\eta \asymp \widehat \eta$  and let $J$ be an arbitrary ideal. Then  $J_{a^2} \supset I_{a^2} \Rightarrow J_a \supset I_a$.
\eT

\noindent \textbf{Question 9.}
Find necessary conditions, at least for the principal ideal case, for the inclusion cancellations (i) and (ii) to hold.\\

By the fact that if $\eta$ satisfies the exponential $\Delta_2$-condition then  $\eta_{a^p}$ satisfies the exponential $\Delta_2$-condition for every $p \in \mathbb N$ \cite[Corollary 3.12] {vKgW04-Density}, Theorems \ref{T: Theorem 4.8},  \ref{T: Theorem 4.5}, and  \ref{T: Theorem 4.6} have a partial extension to higher order arithmetic means.

\begin{theorem}\label{T: higher order cancellation}
Let $p \in \mathbb N$. If every $\co*$-sequence in the characteristic set $\Sigma (I)$ of an ideal $I$ is dominated by some $\eta$ in its characteristic set that satisfies the $p^{th}$ order exponential $\Delta_2$-condition $\underset{m} {\sup}\,\frac{m^2(\eta_{a^p})_{m^2}}{m(\eta_{a^p})_m} < \infty$,
then  $J_{a^{p+1}} = I_{a^{p+1}}$ implies $J_{a^p} = I_{a^p}$, $J_{a^{p+1}} \subset I_{a^{p+1}}$ implies $J_{a^p} \subset I_{a^p}$, and$J_{a^{p+1}} \supset I_{a^{p+1}}$ implies $J_{a^p} \supset I_{a^p}$.
\end{theorem}

The theory of arithmetic mean at infinity ideals, in many instances mirrors that of arithmetic mean ideals. 
We have not investigated second order am-$\infty$ cancellations.\\

\noindent \textbf{Question 10.} What is the corresponding theory for second order am-$\infty$ cancellations?\\

\section{\leftline{\bf The Schur-Horn majorization theorem and am-closed ideals}}\label{S: Majorization}
\subsection{\leftline{\bf Motivations and definitions}}\label{S:7.1}
Given an ideal $I$, $\xi \in \Sigma(I^-)$ if and only if $\xi_a\le\eta_a$ for some $\eta\in  \Sigma(I)$.
The inequality $\xi_a\le\eta_a$ means that $ \sum_{j=1}^n\xi_j = \sum_{j=1}^n\eta_j$ for all $n$ and this relation between (monotone non-increasing) sequences is called \textit{majorization}.

More precisely, in the traditional terminology introduced by Hardy, Littlewood, and P\'olya (see \cite {HLP52}),
\bD{D:maj} Let $\xi, \eta \in (\mathbb R^N)^+$ and let $\xi^*,\eta^*$ denote their monotone rearrangement, then $\xi$ is said to be weakly majorized by $\eta$ if  $ \sum_{j=1}^n(\xi^*)_j \le \sum_{j=1}^n(\eta^*)_j$ for all $1\le n \le N$.
If, in addition, $ \sum_{j=1}^N(\xi^*)_j = \sum_{j=1}^N(\eta^*)_j$, then  $\xi$ is said to be majorized by $\eta$.
\eD
Majorization theory arose during the early part of the 20th century from a number of apparently unrelated  topics: inequalities involving convex functions (Hardy, Littlewood, and P\'olya \cite {HLP52}), wealth distribution (Lorenz \cite {Lm1905}),
convex combinations of permutation matrices (Birkhoff \cite{Bg46}), and more central to our interests \cite{vKgW04-Majorization}, doubly stochastic matrices (Hardy, Littlewood, and P\'olya \cite {HLP52}) and the relation established between eigenvalue lists and diagonals of selfadjoint matrices by the Schur-Horn Theorem:

\begin{theorem}\label{T:S-H} \cite[Schur] {Si23}, \cite[Horn]{Horn}) Let $\xi, \eta \in (\mathbb R^N)^+$. Then  $\xi$ is majorized by $\eta$ if and only if there is a selfadjoint $N\times N$ matrix $A$ having eigenvalue list $\eta$ and diagonal entries $\xi$.
\eT

For infinite sequences ``equality at the end" has no obvious meaning beyond the summable case.  
We provide in the definition below a condition that generalizes this notion for possibly nonsummable infinite sequences which we call \textit{strong majorization}.
We use the following terminology and to simplify notations we introduce it directly for $\co*$ sequences.

\bD{D:majoriz}
Let $\xi, \eta \in \co*$, then we say that
\begin{itemize}
\item
$\xi$ is \textit{majorized} by $\eta$ (denoted by $\xi \prec \eta$) if $ \, \sum_{j=1}^n \xi_j \le \sum_{j=1}^n \eta_j $ for every $n\in \mathbb N; $
\item
 $\xi$ is \textit{strongly majorized} by $\eta$  (denoted by $
 \xi\preccurlyeq \eta)   \text { if } \xi \prec \eta  \, \text { and } {\varliminf  }_{n}\sum_{j=1}^n ( \eta_j -\xi_j) = 0;$
\item $\xi$ is \textit{block majorized} by $\eta$  (denoted by $\xi\prec_b  \eta)  \text { if } \xi \prec \eta  \, \text { and } \sum_{j=1}^{n_k} \xi_j = \sum_{j=1}^{n_k} \eta_j$
for some subsequence $ \mathbb N \ni n_k \uparrow \infty$.
\end{itemize}
\noindent Let $\xi, \eta \in (\ell^1)^*$ then  we say that
\begin{itemize}
\item $\xi$ is \textit{majorized at infinity} by $\eta$  (denoted by $ \xi \prec_\infty\eta$) if $ \sum_{j=n}^\infty \xi_j \le \sum_{j=n}^\infty \eta_j$ for every $ n\in \mathbb N;$
\item  $\xi$ is \textit{strongly majorized at infinity} by $\eta$ (denoted by $ \xi \preccurlyeq_\infty \eta$) if $\xi \prec_\infty\eta $ and \\$ \sum_{j=1}^\infty  \xi_j = \sum_{j=1}^\infty\eta_j$.
\end{itemize}
\eD

Immediate consequences of Definition \ref {D:majoriz} are:
\[
\xi\prec_b  \eta \Rightarrow \xi\preccurlyeq \eta  \Rightarrow \xi \prec \eta,
\]
\[
\text{if }\eta\in(\ell^1)^*\text{ then }\xi\preccurlyeq \eta\, \Leftrightarrow \,\xi\prec \eta \text{ and }\sum_{j=1}^\infty  \xi_j = \sum_{j=1}^\infty\eta_j,
\]
 and
\[
\text{if } \sum_{j=1}^\infty  \xi_j = \sum_{j=1}^\infty\eta_j < \infty, \quad \text{then}\quad  \eta \prec_\infty\xi \Leftrightarrow \xi \prec \eta.
\]
Thus $\xi\in \Sigma (I^-)$ (resp., $\Sigma (I^{-\infty}$) precisely if $\xi\prec\eta$ (resp., $\xi\prec_\infty\eta$) for  some $\eta\in \Sigma(I)$ and $I$ is am-closed (resp., am-$\infty$ closed) precisely when $\Sigma (I)$  is hereditary (i.e., solid) under the majorization order (resp., majorization at infinity order).

Thus majorization at infinity for summable sequences, i.e.,  ``tail majorization,"  that has been introduced and studied also for finite sequences, has a particular relevance for us. Block-majorization can be seen as a way to bring to bear on infinite sequences the results of finite majorization theory, in particular the Schur-Horn Theorem, but we have discovered that it plays also a relevant role on its own.

\subsection{\leftline{\bf Why ``strong" majorization?}}\label{S:7.2}
Beyond the obvious reason that for summable sequences the relation $\preccurlyeq $ is indeed equivalent to majorization with ``equality at the end,"  the relation $\preccurlyeq $ satisfies the following properties that characterize  ``equality at the end" for the finite sequence case and is somewhat reminiscent of density properties discussed earlier in Section \ref{S:5}.

\bT{T:Theorem 5.4}\cite[Theorem 5.4]{vKgW04-Majorization}
Let $\xi, \eta \in \co*$ and $\xi \prec  \eta$.
\item [(i)] There is  a  $\zeta \in \co*$ for which $\xi\preccurlyeq \zeta \le \eta$.
\item [(ii)] There is a  $\rho \in \co*$ for which $\xi \le \rho\preccurlyeq \eta$.
\eT

\bT{T:Theorem 5.7}\cite[Theorem 5.7]{vKgW04-Majorization}
If $\xi,\eta \in (\ell^1)^*$ and   $ \xi \prec_\infty \eta$, then
\item [(i)] $\xi \preccurlyeq_\infty \zeta \le \eta$ for some $\zeta \in \co*$;
\item [(ii)] $\xi \le \rho\preccurlyeq_\infty \eta$ for some $\rho \in \co*$.
\eT

The key tool in the proofs of these results was the following:

\bP {P:Proposition 5.2}\cite [Proposition 5.2] {vKgW04-Majorization}\\
If $\xi, \eta \in \co*$ and $\xi \prec  \eta$, then the following conditions are equivalent.
\item[(i)] For every $m\in \mathbb N$, $\{ \sum_1^n(\eta_j-\xi_j) \mid n \ge m\}$ attains a minimum.
\item[(ii)]   There is  $\zeta \in \co*$ for which $\xi\prec_b  \zeta \le \eta$.
 \item[(iii)] There is  $\rho \in \co*$ for which $\xi \le \rho\prec_b  \eta$.
\eP
\noindent and an analogous but more complicated result helped us handle the "tail condition" \cite [Proposition 5.6]{vKgW04-Majorization}.

 \subsection{\leftline{\bf Infinite majorization and stochastic matrices}}\label{S:7.3}
 Since  Hardy, Littlewood, and P\'olya  \cite {HLP52}, stochastic matrices played a key role in majorization theory. They proved that for finite monotone sequences $\xi$ and $\eta$, $\xi\prec\eta$ (resp., $\xi\preccurlyeq \eta$ -using the notations of our paper, not the original notations of Hardy, Littlewood, and P\'olya) if and only if $\xi=P\eta$ for some substochastic (resp., doubly stochastic) matrix $P$. These are matrices with non-negative entries, where all the rows and columns have sums $\le 1$ (resp., $=1$). Closer to our interest, the key step in the proof of the Schur-Horn theorem was to show that if $\xi\preccurlyeq \eta$ then $\xi=Q\eta$ for some orthostochastic matrix $Q$, i.e., a matrix $Q_{ij}=(U_{ij})^2$ for some orthogonal matrix $U$ (unitary with real entries). Notice that it is then an elementary computation to show that
 $\diag \xi =E(U\diag \eta\, U^*)$ where $E(X)$ denotes the main diagonal of a matrix $X$.

Majorization for infinite sequences has received over the years only a moderate amount of attention.
In a not nearly as well-known article as it deserves, in 1964 Markus \cite [Lemma 3.1]{aM64}  extended to the infinite case part of the Hardy, Littlewood, and P\'olya theorem, by  showing that if $\xi, \eta \in \co*$, then $\xi \prec \eta $ if and only if $\xi = Q\eta$ for some substochastic matrix $Q$. Although he did not adopt this terminology, his proof was based on iterations of T-transforms (see \cite[Section 4]{vKgW04-Majorization} for details).

In fact, with a slight tightening of his construction, it is easy to see that when  $\xi_n>0$ for all $n$, then the  construction associates canonically with $\xi \prec \eta $ a substochastic matrix $Q(\xi, \eta)$ for which $\xi = Q(\xi, \eta)\eta$ and as Markus had already remarked, the matrix $Q(\xi, \eta)$ is row-stochastic (all the rows have sum equal to 1). What was, however, much less trivial to show was

\bT{T: Theorem 3.4}\cite [Theorem 3.4]{vKgW04-Majorization}  Let $\xi, \eta \in \co*$ with $\xi_n > 0$ for every $n$ and  $\xi \prec \eta$. 
Then  $Q(\xi,\eta)$ is co-isometry stochastic, i.e., $Q(\xi,\eta)_{ij}=(W(\xi,\eta)_{ij})^2$ and $W(\xi,\eta)$ is a co-isometry  with real entries.
\eT

The co-isometry $W(\xi,\eta)$ too is canonical and completely determined by a sequence $\{m_k,t_k\}$ where $m_k\in \mathbb N$ are the size of the T-transform matrices and $0 < t_k \le 1$ the convex coefficients (see \cite{vKgW04-Majorization} for details). We further have:

\bT{T:Theorem 3.3}\cite[Theorem 3.3]{vKgW04-Majorization}
 Let  $\xi \prec \eta $  for some $\xi, \eta \in \co*$ with $\xi_n >0$ for every $n$. Then the following are equivalent.
\item[(i)] $\xi\preccurlyeq   \eta$
\item[(ii)] $W(\xi,\eta)$ is an orthogonal matrix.
\item[(iii)]  $\sum\{t_k \mid  m_k=1\}=\infty$.
\eT

\bC{C:Corollary 3.13}\cite[Corollary 3.13]{vKgW04-Majorization} Let $\xi \prec \eta $  for some $\xi, \eta \in \co*$ with $\xi_n >0$ for every $n$. 
Then the following are equivalent.
\item[(i)] $\xi\prec_b  \eta$
\item[(ii)] $W(\xi,\eta)$ is the direct sum of finite orthogonal matrices.
\item[(iii)] The set $\{k\mid m_k=t_k=1\}$ is infinite.
\eC

Theorem \ref {T:Theorem 3.3} shows that if $\xi\preccurlyeq \eta$ then  $\xi=Q\eta$ for some orthostochastic matrix $Q$, namely $Q=Q(\xi,\eta)$. We found an example of an nonsummable $\xi$ with $\xi\prec \eta$ but $\xi\not \preccurlyeq \eta$ where  $\xi=Q\eta$ for some $Q\ne Q(\xi,\eta)$ orthostochastic \cite[Example 2.12 ]{vKgW04-Majorization}.  It turned out that this example was ``generic" (see Theorem \ref {T: Theorem 4.3} below). The key step in the proof was  the following ``splice" lemma:

\bL{L:Lemma 4.2} \cite[Lemma 4.2]{vKgW04-Majorization} Let $\xi,\eta\in \co*$ and assume that  $\xi\prec \eta $ but $\xi \not\preccurlyeq \eta$ and $\xi \not\in (\ell^1)^*$. Then there is a partition of $\mathbb N$ in two sequences $\{n_j^{(1)}\}$ and $\{n_j^{(2)}\}$ with $n_1^{(1)}=1$, so that if $\xi':=\{\xi_{n_j^{(1)}}\}$, $\eta':=\{\eta_{n_j^{(1)}}\}$, $\xi'':=\{\xi_{n_j^{(2)}}\}$, $\eta'':=\{\eta_{n_j^{(2)}}\}$ are the corresponding subsequences of $\xi$ and $\eta$, then $\xi'\preccurlyeq \eta'$ and $\xi''\prec \eta''$ but $\xi''\not\preccurlyeq \eta''$ and $\xi''\not\in (\ell^1)^*$.
\eL
Then a direct sum argument yields

\bT{T: Theorem 4.3}\cite[Theorem 4.3]{vKgW04-Majorization}\\
If $\xi, \eta\in \co*$ and $\xi\not\in (\ell^1)^*$, then $\xi\prec\eta$ if and only if $\xi=Q\eta$ for some  orthostochastic matrix $Q$.
\eT
\noindent which in turns is the key for our infinite dimensional Schur-Horn majorization theorem:

 \bT{T:New Horn}\cite[Theorem 4.4]{vKgW04-Majorization}\\
Let $\xi, \eta\in \co*$. Then the following conditions are equivalent.
\item[(i)] There is an $A\in K(H)^+$ with $s(A)=\eta$ and $s(E(A))=\xi$.
\item[(ii)] $\diag \xi = E(U\diag\eta U^*)$ for some unitary operator $U$.
\item[(iii)] $\xi=Q\eta$ for some  orthostochastic matrix $Q$. \\
\noindent If $\xi\in (\ell^1)^*$, then conditions (i)-(iii) are equivalent to
\item[(ivS)]  $\xi\preccurlyeq \eta$\\
\noindent If $\xi\not\in (\ell^1)^*$, then conditions (i)-(iii) are equivalent to
\item[(ivNS)]  $\xi\prec \eta$.
\eT

The nonsummable case is new. The summable case is an extension of the similar result by Arveson and Kadison \cite[Theorem 4.1]{AK02}.  

There has been a resurgence of interest in various aspects of the Schur-Horn Theorem and in particular its convexity formulation in terms of groups.  The important paper by Konstant \cite {Kb73} was extended to infinite dimensions by Block, Flaschka, and Ratiu \cite{BFR93} and by A. Neumann  \cite{Na99}. Neumann did also extend in part the connection between majorization theory (but for more general sequences) and diagonal of matrices, and further work in this direction was obtained  by Antezana,  Massey,  Ruiz, and Stojanoff \cite{AMRS}. 

There is fairly little overlap between our results and those in \cite{Na99}  and \cite{AMRS}. 
While we focus on individual operators with prescribed eigenvalue list and diagonal, the cited papers deal with norm closures of classes of such operators.
Some connections and differences are explored in \cite[Section 6]{vKgW04-Majorization}.

\subsection{\leftline {Applications to operator ideals}} \label{S:6.1} \quad
The properties of majorization for infinite sequences developed above permit to obtain the following characterization of the am-closure (resp., am-$\infty$ closure) of an ideal.
\bT{T:Theorem 6.1}\cite [Theorem 6.1]{vKgW04-Majorization}\\
Let $I$ be an ideal and let $ \xi \in \co*$. Then the following conditions are equivalent.
\item[(i)] $\xi \in \Sigma (I^-)$.
\item[(ii)]  $\xi \prec \eta$ for some $\eta \in  \Sigma(I)$.
\item[(ii$'$)] $\xi \preccurlyeq \eta$ for some $\eta \in  \Sigma(I)$.
\item[(iii)] $\xi =P \eta$ for some $\eta \in  \Sigma(I)$ and some substochastic matrix $P$.
\item[(iii$'$)]  $\xi =Q \eta$ for some $\eta \in  \Sigma(I)$ and some orthostochastic matrix $Q$.
\item[(iv)]   $\diag \xi = E(L \diag  \eta\, L^*)$ for some $\eta \in \Sigma(I)$ and some contraction $L$.
\item[(iv$'$)]  $\diag \xi = E(U \diag  \eta\, U^*)$ for some $\eta \in  \Sigma(I)$ and some orthogonal matrix $U$.\\
If $I\supset \mathscr L_1$, then the above conditions are further equivalent to
\item[(ii$''$)] $\xi\prec_b  \eta$ for some $\eta \in  \Sigma(I)$.
\item[(iii$'$)]  $\xi =Q \eta$ for some $\eta \in  \Sigma(I)$ and some block orthostochastic matrix $Q$.
\item[(iv$''$)]  $\diag \xi = E(U \diag  \eta\, U^*)$ for some $\eta \in  \Sigma(I)$ and some matrix $U$ direct sum of finite orthogonal matrices.
\eT
\bT{T:Theorem 6.6}\cite [Theorem 6.6]{vKgW04-Majorization}
Let $I\subset \mathscr L_1$ be an ideal and let $ \xi \in (\ell^1)^*$. Then the following conditions are equivalent
\item[(i)] $\xi \in \Sigma (I^{-\infty})$
\item[(ii)]  $ \xi \prec_\infty\eta$ for some $\eta \in  \Sigma(I)$.
\item[(ii$'$)]  $\xi \preccurlyeq_\infty \eta$ for some $\eta \in  \Sigma(I)$.
\item[(ii$''$)] $\eta \preccurlyeq \xi$ for some $\eta \in  \Sigma(I)$.
\item[(iii)]  $P\xi \in \Sigma(I)$ for some column stochastic matrix $P$.
\item[(iii$'$)]  $Q\xi \in \Sigma(I)$ for some block orthostochastic matrix $Q$.
\item[(iv)]   $E(V \diag  \xi\, V^*)\in I$ for some isometry $V$.
\item[(iv$'$)]   $E(U \diag  \xi\, U^*)\in I$ for some matrix $U$ direct sum of finite orthogonal matrices.
\eT
As a consequence, if we let $\mathscr D$ denote the masa of diagonal operators (for the given orthonormal basis), we have:

\bC{C:amclosure}
For every ideal $I$, $E(I) = I^-\cap \mathscr D$.
\eC
and hence
\bC{C:amclosed}
An ideal $I$ is am-closed if and only if $E(I) \subset I$.
\eC
The latter property was called diagonal invariance in \cite {gW75} and it can be viewed as a generalization of \cite [Theorem 4.2, Chapter III] {GK69} which proves that $E(\mathfrak{S}_\phi)\subset \mathfrak{S}_\phi$ for every maximal Banach ideal  $\mathfrak{S}_\phi$, where $\phi$ is the generating 
symmetric norming function.

Theorem \ref {T:Theorem 6.1} provides a further characterization of am-closed ideals:

\bC{C:Corollary 6.7} Let $ I$ be an ideal. Then the following are equivalent:
\item[(i)] $I$ is am-closed.
\item[(ii)] $\Sigma(I)$ is invariant under substochastic matrices, in the sense that if $\xi\in \Sigma(I)$ and $P$ is a substochastic matrix, then $(P\xi)^*\in \Sigma(I)$, where * denotes  monotone rearrangement.
\item[(iii)] $\Sigma(I)$ is invariant under orthostocastic matrices and hence under any of the following classes of (sub)stochastic matrices: row-stochastic, column-stochastic, doubly stochastic, isometry stochastic, co-isometry stochastic, unitary  stochastic.\\
\noindent If $I\supset \mathscr L_1$, then the above conditions are also equivalent to:
\item[(iv)] $\Sigma(I)$ is invariant under block stochastic matrices.
\eC

Without the condition $I\supset \mathscr L_1$, (iv) does not  imply (i), e.g., $F$ is invariant under block stochastic matrices but it is not am-closed. More generally, in \cite{vKgW04-Soft} we have shown that every am-closed ideal must contain  $\mathscr L_1$.

Invariance under various classes of substochastic matrices can be seen as an (infinite) convexity property  of the  ideal and the following will make this clearer.
Kendall introduced in \cite{Kd60} a topology on the space of boundedly row and column summable infinite real matrices with respect to which the set of doubly stochastic matrices is precisely the closed convex hull of the set of permutation matrices and used it to prove that the extreme points of the set of doubly stochastic matrices are the permutation matrices.

The class of infinite convex combinations of permutation matrices
\[
\mathscr C:=\{\sum_{j=1}^\infty t_j\Pi_j \mid \Pi_j \text { permutation matrix, } 0\le t_j\le1,~  \sum_{j=1}^\infty t_j=1\}
\]
is strictly smaller than the class of doubly stochastic matrices. This  can be seen for instance by noticing that $\mathscr C\cap K(H)=\{0\}$ while it is possible to find a Hilbert-Schmidt doubly stochastic matrix \cite{mCgW02}. More interestingly, invariance of a characteristic set $\Sigma(I)$ of an ideal  ideal $I$ under $\mathscr C$ does not imply invariance of $\Sigma(I)$ under the class of doubly stochastic matrices, i.e., am-closure of $I$. Indeed, it elementary to see that If $I$ is a Banach ideal then $\Sigma(I)$ is invariant under $\mathscr C$ and yet there are Banach ideals that are not am-closed, for example the closure $\cl(\xi)$ of a principal ideal $\xi$ for $\xi$ irregular and nonsummable 
\cite {jV89} (see also discussion after Proposition \ref{P:Gohberg/Krein soft pairs}). It is worth noticing that this counterexample is not a soft ideal, and this seemingly unrelated fact turns out to be quite relevant as we see in the next proposition,

Let  class $\mathscr B$ denote the following class of block substochastic  matrices:
\[
\mathscr B:=\{\sum_{j=1}^\infty \oplus t_j\Pi_j \mid \Pi_j \text { finite permutation matrix, } 0\le t_j\le1,~  \sum_{j=1}^\infty t_j=1\}.
\]

\bP{P:(i),(ii))} \cite [Theorem 6.9] {vKgW04-Majorization} Let $I$ be an ideal.
\item[(i)] If $\Sigma(I)$ is invariant under $\mathscr C$, then $I\supset \mathscr L_1$ and $\Sigma(I)$ is invariant under $\mathscr B$.
\item[(ii)] If $I$ is soft-edged or soft-complemented, and $\Sigma(I)$ is invariant under $\mathscr B$, then $I$ is invariant under block stochastic matrices. 
\item[(iii)] If $I$ is soft-edged or soft-complemented and $\Sigma(I)$ is invariant under $\mathscr C$, then $I$ is am-closed.
\eP


\begin{thebibliography}{99}

\bibitem{AGPS87}
Albeverio, S., Guido, D., Ponosov, A., and Scarlatti, S., 
\textit{Singular traces and compact operators,} J. Funct. Anal. \textbf{137}(2) (1996), pp.~281--302.

\bibitem{AS78}
G. D. Allen and L. C. Shen, \textit{On the structure of principal
ideals of operators.} Trans. Amer. Math. Soc. \textbf{238} (1978),
253--270.


\bibitem{jA77}
Anderson, J., \textit{Commutators of compact operators,} J. Reine Angew. Math. \textbf{291} (1977), pp.~128--132.

\bibitem{AV86}
Anderson, J. and Vaserstein, L. N., \textit{Commutators in ideals of trace class operators,} Indiana University Mathematics Journal \textbf{35} (2) (1986), pp.~345--372.

\bibitem{jA86}
Anderson, J., \textit{Commutators in ideals of trace class operators. {II},} Indiana University Mathematics Journal \textbf{35} (2) (1986), pp.~373--378.


\bibitem{AMRS}
Antezana,  J.,  Massey, P.,  Ruiz, M, and Stojanoff, D., \textit{The Schur-Horn Theorem for operators and frames with prescribed norms and frame operator,} Preprint.

\bibitem{AK02}
Arveson, W., and Kadison, R. V.,  \textit{Diagonals of self-adjoint operators.} Preprint.

\bibitem{Bg46}
Birkhoff, G., \textit{Tres observaciones sobre el algebra lineal,} Univ. Nac. Tucuman Rev. Ser. A (\textbf{5}) (1946) 147-151

\bibitem{aB89}
Blass, A., \textit{Applications of superperfect forcing and its
relatives,} Set Theory and its Applications, Lecture Notes in
Mathematics 1401, Springer-Verlag (1989), 18-40.

\bibitem{aB90}
Blass, A., \textit{Groupwise density and related cardinals,} Arch. Math. Logic 30 (1990) 1-11.

\bibitem{aBgW78}
Blass, A., and Weiss, G.,  \textit{A Characterization and Sum
Decomposition of Operator Ideals,} Trans. Amer. Math. Soc.
\textbf{246} (1978), 407-417.

\bibitem{BFR93} 
Block, A. M.,  Flaschka, H.,  and Ratiu, T. S., \textit{A Schur-Horn-Konstant convexity theorem for the diffeomorphism group of the annulus,} Invent. Math. \textbf{113} no 3, (1993), 511-529.


\bibitem{BP65}
Brown, A. and Pearcy, C., \textit{Structure of commutators of operators,} Ann. of Math. (2) \textbf{82} (1965), pp.~112--127.

\bibitem{BPS71} Brown, A., Pearcy, C., and Salinas, N.,
\textit{Ideals of compact operators on {H}ilbert space,} Trans. Amer.
Math. Soc. \textbf{18} (1971), 373--384.

\bibitem{jC41}
Calkin, J. W., \textit{Two-sided ideals and congruences in the ring
of bounded operators in {H}ilbert space,}
Ann. of Math. (2) \textbf{42} (1941), pp.~839--873.

\bibitem{mCgW02}
Chance, M. and Weiss, G., private communication, 2002.

\bibitem{aC94}
Connes, A., \textit{Non Commutative Geometry,} San Diego Academic Press, 1994.

\bibitem{aC82}
Connes, A., \textit{Noncommutative differential geometry Part I, the Chern character in $K$-homology,} 
Inst. Hautes Etudes Sci., Bures-Sur-Yvette, 1982.

\bibitem{aC83}
Connes, A., \textit{Noncommutative differential geometry Part II, de Rham homology and noncommutative algebra,} 
Inst. Hautes Etudes Sci., Bures-Sur-Yvette, 1983.

\bibitem{aC85}
Connes, A., \textit{Noncommutative differential geometry,} 
Inst. Hautes \'Etudes Sci. Publ. Math. \textbf{62} (1985), pp.~257--360. 

\bibitem{jD57}
Dixmier, J., \textit{Les alg\`ebres d'op\'erateurs dans l'espace hilbertien ({A}lg\`ebres de von {N}eumann)}, 
Cahiers scientifiques, fasc. 25, Gauthier-Villars, Paris, 1957. 


\bibitem{jD66}
Dixmier, J., \textit{Existence de traces non normales,} 
C. R. Acad. Sci. Paris S\'er. A-B \textbf{262} (1966), pp.~A1107--A1108.


\bibitem{kDgWmW00}
Dykema, K., Weiss, G., and Wodzicki, M., \textit{Unitarily invariant trace extensions beyond the trace class,} 
Complex analysis and related topics (Cuernavaca, 1996), Oper. Theory Adv. Appl. \textbf{114} (2000), pp.~59--65.


\bibitem{DFWW}
Dykema, K., Figiel, T., Weiss, G., and Wodzicki, M., \textit{The
commutator structure of operator ideals,}
Adv. Math., 185/1 pp. 1--79.


\bibitem{kDnK98}
Dykema, K. J. and Kalton, N. J., \textit{Spectral characterization of sums of commutators. {I}{I},} J. Reine Angew. Math. \textbf{504} (1998), pp.~127--137.

\bibitem{kDnK05}
Dykema, K. J. and Kalton, N. J., \textit{Sums of commutators in 
ideals and modules of type {II} factors}
Ann. Inst. Fourier (Grenoble) \textbf{55} 3 (2005), pp.~931--971.

\bibitem{FKM06}
Feldman, I., Krupnik, N., Markus, A. \textit{On the connection between the indices of a block operator matrix and of its determinant,} 
Operator Theory: Advances and Applications \textbf{1}, Birka\"user Verlag Basel/Switzerland (2006), pp.~1--17.


\bibitem{bF51}
Fuglede, B., \textit{A commutativity theorem for normal operators,} Proc. Nat. Acad. Sci. U. S. A. \textbf{36} (1950), pp.~35--40.

\bibitem{GK69}
Gohberg, I. C. and Kre{\u\i}n, M. G., \textit{Introduction to the
theory of linear nonselfadjoint operators.} American Mathematical
Society (1969).


\bibitem{pH54}
Halmos, P. R., \textit{Commutators of operators. {II},} Amer. J. Math. \textbf{76} (1954), pp.~191--198.

\bibitem{HLP52}
Hardy, G. H., Littlewood, J. E. and P\'olya, G., \textit{Inequalities.} 2d ed. Cambridge University Press, 1952.

\bibitem{Horn}
Horn, A. \textit{Doubly stochastic matrices and the diagonal of a rotation matrix.} Amer. J. Math.  \textbf{76} (1954), 620-630.


\bibitem{vKgW02}
Kaftal, V. and Weiss, G., \textit{Traces, ideals, and arithmetic
means,} Proc. Natl. Acad. Sci. USA \textbf{99} (11) (2002),
pp.~7356--7360.

\bibitem{vKgW04-Traces}
Kaftal, V. and Weiss, G., \textit{Traces on operator ideals and
arithmetic means,} preprint.

\bibitem{vKgW04-Soft}
Kaftal, V. and Weiss, G., \textit{Soft ideals and arithmetic mean
ideals,} IEOT, to appear.

\bibitem{vKgW04-Density}
Kaftal, V. and Weiss, G., \textit{The $B(H)$ lattices, density and 
arithmetic mean ideals,} preprint.

\bibitem{vKgW04-2nd Order}
Kaftal, V. and Weiss, G., \textit{Second Order Arithmetic Means in Operator Ideals,} J. Operators and Matrices (OAM), to appear.

\bibitem{vKgW04-Majorization}
Kaftal, V. and Weiss, G., \textit{Majorization for infinite sequences, an extension of the Schur-Horn Theorem, and operator ideals,} preprint.

\bibitem{nK87}
Kalton, N. J., \textit{Unusual traces on operator ideals,} Math. Nachr. \textbf{134} (1987), pp.~119--130. 


\bibitem{nK89}
Kalton, N. J., \textit{Trace-class operators and commutators,} J.
Funct. Anal. \textbf{86} (1989), pp.~41--74.

\bibitem{nK98}
Kalton, N. J., \textit{Spectral characterization of sums of commutators {I},} J. Reine Angew. Math. \textbf{504} (1998),
pp.~115--125.

\bibitem{Kd60}
Kendall, D., \textit{ On infinite doubly-stochastic matrices and Birkhoff's problem 111}, Journal London Math. Soc, \textbf {35} (1960), 81-84.

\bibitem{Kb73}
Konstant, B.,  \textit{On convexity, the Weyl group  and the Iwasawa decomposition,} Ann. Sci. \'{E}cole Norm. Sup. \textbf {4} 6 (1973), 413-455.


\bibitem{Lm1905}
Lorenz, M. O.,  \textit{Methods of measuring concentration of wealth,} J. Amer. Statist. Assoc, \textbf {9} (1905), 209-219. 

\bibitem{aM64}
A. S. Markus, \textit{The eigen- and singular values of the sum and
product of linear operators.} Uspekhi Mat. Nauk \textbf{4} (1964),
93--123.


\bibitem{MO79}
Marshall, A. W. and Olkin, I., \textit{Inequalities: Theory of 
Majorization and its Applications,} Academic Press Inc. [Harcourt 
Brace Jovanovich Publishers], Mathematics in Science and Engineering 
\textbf{143} (1979).

\bibitem{hM01}
Mildenberger, H.,
\textit{Groupwise dense families,} Arch. Math. Logic \textbf{40}
(2001), no. 2, pp. 93--112.

\bibitem{Na99}
Neumann, A.,  \textit{An infinite-dimensional generalization of the Schur-Horn convexity theorem.} Jour. Funct. Anal., \textbf{161} (2), {1999},  418-451.

\bibitem{Oa52}
Ostrowski, A. M. , \textit{Sur quelques applications des fonctions convexes et concaves au sens the I. Schur} J. Math. Pures Appl. [9], \textbf{31}, (1952), 253-292.

\bibitem{PT71}
Pearcy, C. and Topping, D., \textit{On commutators in ideals of compact operators,} Michigan Math. J. \textbf{18} (1971), pp.~247--252.

\bibitem{cP51}
Putnam, C. R., \textit{On normal operators in {H}ilbert space,} Amer. J. Math. \textbf{73} (1951), pp.~357--362.


\bibitem{nS74}
N. Salinas, \textit{Symmetric norm ideals and relative conjugate
ideals.} Trans. Amer. Math. Soc. \textbf{138} (1974), 213--240.

\bibitem{rS60}
Schatten, R., \textit{Norm ideals of completely continuous operators,} Ergebnisse der Mathematik und ihrer Grenzgebiete, 	Neue Folge, Heft~ {\bf 27} , Berlin, Springer-Verlag, 1960.

\bibitem{Si23}
Schur, I., \textit{\"{U}ber eine Klasse von Mittlebildungen mit Anwendungen auf der Determinantentheorie,} Sitzungsber. Berliner Mat. Ges., ( \textbf{22}), (1923),  9-29.

\bibitem{vS83}
Shul'man, V. S., \textit{Linear equations with normal coefficients,} Dokl. Akad. Nauk SSSR \textbf{270} {(5)} (1983), pp.~1070--1073.


\bibitem{jV89}
Varga, J., \textit{Traces on irregular ideals,} Proc. Amer. Math.
Soc. \textbf{107} 3 (1989), pp.~715--723.

\bibitem{dV79/81}
Voiculescu, D., \textit{Some results on norm-ideal perturbations of {H}ilbert space operators {I}{I},} 
J. Operator Theory \textbf{5} (1981), pp.~77--100.


\bibitem{gW75}
Weiss, G., \textit{Commutators and Operators Ideals,} dissertation
(1975), University of Michigan Microfilm.


\bibitem{gW80}
Weiss, G., \textit{Commutators of {H}ilbert-{S}chmidt operators.
{I}{I},} IEOT \textbf{3} (4) (1980), pp.~574--600.

\bibitem{gW83}
Weiss, G., \textit{An extension of the {F}uglede commutativity theorem modulo the
              {H}ilbert-{S}chmidt class to operators of the form $\sum
              {M}\sb{n}{X}{N}\sb{n}$,} Trans. Amer. Math. Soc. \textbf{278} (1) (1983), pp.~1--20.

\bibitem{gW86}
Weiss, G., \textit{Commutators of {H}ilbert-{S}chmidt operators.
{I},} IEOT \textbf{9} (6) (1986), pp.~877--892.

\bibitem{gW05}
Weiss, G., \textit{B(H)-commutators: a historical survey.} Recent advances in operator theory, operator algebras, and their applications, 307--320, Oper. Theory Adv. Appl., 153, Birkhauser, Basel, 2005.


\bibitem{mW94}
Wodzicki, M., \textit{Algebraic {$K$}-theory and functional analysis,} 
First European Congress of Mathematics, Vol.\ II (Paris,1992) \textbf{120} (1994), pp.~485--496. 

 
\bibitem{mW02}
Wodzicki, M., \textit{Vestigia investiganda,} Mosc. Math. J. \textbf{4} (2002), pp.~769--798. 

\end{thebibliography}
\end{document}